\documentclass[11pt,twoside]{article}
\usepackage{latexsym}
\usepackage{amssymb,amsbsy,amsmath,amsfonts,amssymb,amscd}
\usepackage{mathtools}
\usepackage{graphicx, color}
\usepackage[colorlinks=true]{hyperref}

\newcommand{\fer}[1]{(\ref{#1})}

\usepackage{subfig}
\usepackage{caption}
\newcommand{\commentout}[1]{}
\newcommand{\R}{\mathbb{R}}

\newcommand {\e}  {\varepsilon}

\newcommand {\f}   {\frac}
\newcommand {\p}   {\partial}
\newcommand{\dis}{\displaystyle}
\newcommand {\proof} {\noindent {\bf Proof}. }
\newcommand{\beq}{\begin{equation}}
\newcommand{\eeq}{\end{equation}}
\newcommand{\bea} {\begin{array}{rl}}
\newcommand{\eea} {\end{array}}
\newcommand{\bepa}{\left\{ \begin{array}{l}}
\newcommand{\eepa} {\end{array}\right.}
\newcommand{\na}{\nabla}
\newcommand{\ep}{\epsilon}
\newcommand{\ddt}{\frac{d}{dt}}
\newcommand{\re}{\eqref}
\newcommand{\De}{\Delta}
\newtheorem{theorem}{Theorem}[section]
\newtheorem{lemma}[theorem]{Lemma}

\newcommand{\qed}{{ \hfill
                    {\unskip\kern 6pt\penalty 500 \raise -2pt\hbox{\vrule\vbox to 6pt{\hrule width 6pt
                    \vfill\hrule}\vrule} \par}   }}
\title{\Large \bf Dirac mass dynamics in multidimensional nonlocal parabolic equations}

\author{
Alexander Lorz\thanks{Department of Applied Mathematics and Theoretical Physics (DAMTP),
Centre for Mathematical Sciences, Wilberforce Road,
Cambridge CB3 0WA, UK. Email: A.Lorz@damtp.cam.ac.uk}
\and
Sepideh Mirrahimi
\thanks{UPMC, CNRS UMR 7598, Laboratoire Jacques-Louis Lions, F-75005, Paris. Email: mirrahimi@ann.jussieu.fr}
\and Beno\^ \i t Perthame \footnotemark[2] \thanks{and Institut Universitaire de France. Email: benoit.perthame@upmc.fr}
}

\usepackage{fancyhdr}

\numberwithin{equation}{section}

\begin{document}
\maketitle


\pagenumbering{arabic}

\begin{abstract}
Nonlocal Lotka-Volterra models have the property that solutions concentrate as Dirac masses in the limit of small diffusion.  Is it possible to describe the dynamics of the limiting concentration points and of the weights of the Dirac masses? What is the long time asymptotics of these Dirac masses? Can several Dirac masses co-exist?

We will explain how these questions relate to the so-called ''constrained Hamilton-Jacobi equation'' and how a form of canonical equation can be established.  This equation has been established assuming smoothness. Here we build a framework where smooth solutions exist and thus the full theory can be developed rigorously. We also show that our form of canonical equation comes with a kind of Lyapunov functional.

Numerical simulations show that the trajectories can exhibit unexpected dynamics well explained by this equation.

Our motivation comes from population adaptive evolution a branch of mathematical ecology which models darwinian evolution.\\

This is a new version of the article published in CPDE in 2011. We have included an additional assumption~\eqref{as:survival} which is  needed for the  Theorem~\ref{th.conv}. We use this assumption  to provide a lower bound for the total population size in Section~\ref{se.rho}. Once this assumption is made, all the previous arguments used in the  original article hold. One of the authors  realized this subtlety while preparing a recent article~\cite{MC.CE.SM:19} which provides precise conditions for the extinction or survival of the population.

\end{abstract}

\section{Motivation}
\label{sec:motivation}

The nonlocal Lotka-Volterra parabolic equations arise in several areas such as ecology,  adaptative dynamics and can be 
derived from stochastic individual based models in the limit of infinite population. The simplest example assumes competition
 between individuals with a trait $x$, through a single resource and reads
\begin{equation} \label{para}
\p_t n_\ep - \ep \Delta n_\ep=  \frac{n_\ep}{\ep} R \big(x, I_\ep(t) \big),  \qquad t>0, \; x \in \R^d,
\end{equation}
with a nonlinearity driven by the integral term
\begin{equation} \label{paraI}
I_\ep(t) = \int_{\R^d}  \psi(x) n_\ep(t,x) dx.
\end{equation}
Another and more interesting example is with direct competition
\begin{align}\label{localc}
\p_t n_\ep(t,x) = \frac{1}{\ep}n_\ep(t,x) \left(r(x)- \int_{\R^d} C(x,y)n_\ep(t,y) \,dy\right) + \ep \De n_\ep(t,x).
\end{align}
We denote by $n_\ep^0\geq 0$ the initial data.

These are called 'mutation-competition' models because the Laplace term is used for modeling mutations in the population. Competition is taken into account in the second model by the competition kernel $C(x,y)\geq 0$ and in the first model by saying that $R$ can be negative for $I_\ep$ large enough (it is a measure of how the total population influences birth and death rates). Such models can be derived from stochastic individual based models in the limit of large populations, \cite{RC.JH:05, NC.RF.SM:06,NC.RF.SM:08}.  There is a large literature on the subject, in terms of modeling and analysis, we just refer the interested reader to \cite{OD:04,OD.PJ.SM.BP:05,M.P.B.M:10,Raoulphd}.
\\

We have already normalized the model with a small positive parameter $\e$ since it is our goal to study the behaviour of the solution as $\e \to 0$. The interesting qualitative outcome is that solutions concentrate as Dirac masses
$$
n_\ep (t,x) \approx \bar \rho(t) \delta \big(x- \bar x (t)\big).
$$
For equation \eqref{para}, we can give an intuitive explanation; in this limit, we expect that the relation $n(t,x) R \big(x, I(t) \big)=0$ holds. In dimension 1 and for $x \mapsto R(x,I)$ monotonic, there is a single point $x=X(I)$ where $R$ will vanish and, consequently, where $n$ will not vanish. A priori control of the total mass on $n$ from below implies the result with $\bar x (t)= X(I(t))$.

In several studies, we have established these singular limits with weak assumptions \cite{GB.SM.BP:09, GB.BP:08}. A main new concept arises in this limit, the {\em constrained Hamilton-Jacobi} equation introduced in \cite{OD.PJ.SM.BP:05}  which occurs by some kind of real phase WKB ansatz (as for fronts propagations in \cite{WF.PS:86, LE.PS:89,GB.LE.PS:90})
\begin{equation} \label{wkb}
n_\ep (t,x) = e^{u_\ep(t,x)/\ep}.
\end{equation}
Here we have in mind the simple example of Dirac masses approximated by gaussians
$$
\delta(x-\bar x) \approx \f {1}{\sqrt{2\pi \ep}} e^{-|x-\bar x|^2/2\ep}= e^{(-|x-\bar x|^2- \ep ln(2 \pi\ep))/2\ep}.
$$
It is much easier to describe the limit of $-|x-\bar x|^2- \ep ln(2 \pi\ep)$! Dirac concentration points are understood as maximum points of $u_\ep(t,x)$ in \eqref{wkb}.
As it is well understood, these Hamilton-Jacobi equations develop singularities in finite time \cite{GB:94,LE:98,Lionsbook:82} which is a major technical difficulty both for proving the limit and for analyzing properties of the concentration points $\bar x (t)$.
\\

This method in \cite{OD.PJ.SM.BP:05}  of using $\ep \ln (n_\ep)$ to prove concentration has been followed in several subsequent studies. For long time asymptotics (and not $\ep \to 0$ but the two issues are connected as we explain in section \ref{sec:senm}) it was introduced in \cite{LD.PJ.SM.GR:08} and used in \cite{Raoul2009,Raoulphd}. More recently in \cite{NC.PJ:10} the authors come back on the Hamilton-Jacobi equation and prove that it makes sense still with weak assumptions for several nonlocal quantities $I_k=\int \psi_k n_\ep(t,x)dx$ which can be characterized in the limit.
\\

Here we take the counterpart and develop a framework where we can prove smoothness of the various quantities arising in the theory.  This opens up the possibility to address many questions that seem impossible to attain directly
\\
$\bullet$ Do the Dirac concentrations points appear spontaneously at their optimal location or do they move regularly?
\\
$\bullet$ In the later case, is there a differential equation  on the concentration point $\bar x(t)$?  It follows from regularity that we can establish a form of the so-called {\em canonical equation} in the language of adaptive dynamics \cite{OD:04,UD.RL:96}. This equation has been established {\em assuming} smoothness in \cite{OD.PJ.SM.BP:05}, in our framework it holds true.
\\
$\bullet$ In higher dimensions, why is a single Dirac mass naturally sustained (and not the hypersurface $R(\cdot,I)=0$ for instance)? The canonical equation enforces constraints on the dynamics which give the explanation.
\\
$\bullet$ What is the long time behaviour of the concentration points $\bar x(t)$? A simple route is that the canonical equation comes with a kind of Lyapunov functional.
\\

We develop the theory separately for the simpler case of the model with competition through a single resource \eqref{para} and for the direct competition
 model \eqref{localc}. For the model with a single resource  we rely on assumptions stated in section \ref{sec:gca} and we give all the details of the proofs 
in the three subsequent sections. We illustrate the results with numerical simulations that are presented in section \ref{sec:numerics}.
 We give several extensions afterwards; in section \ref{sec:extension} we treat the case with non-constant diffusion, and 
finally the case of direct competition in section \ref{sec:localcompetition}.

\section{A simple example: no mutations}
\label{sec:senm}

The Laplace term in the asymptotic analysis of \eqref{para} and \eqref{localc} is at the origin of several assumptions and  technicalities. In order to explain our analysis in a simpler framework, we begin with the  case of the two equations without mutations set for  $t>0$,  $x \in \R^d$,
\begin{equation} \label{paranm}
\p_t n=  {n} R \big(x, I(t) \big),   \qquad I(t)= \int_{\R^d}  \psi(x) n(t,x) dx.
\end{equation}
\begin{align}\label{localcnm}
\p_t n(t,x) = n(t,x) \left(r(x)- \int_{\R^d} C(x,y)n(t,y) \,dy\right):= n(t,x) R \big(x, I(t,x) \big).
\end{align}
Also, we give a formal analysis, that shows the main ideas and avoids writing a list of assumptions; those of the section \ref{sec:gca} and \ref{sec:localcompetition} are enough for our purpose.

In both cases one can easily see the situation of interest for us. The models admit a continuous family of singular, Dirac masses, steady states parametrized by $y\in \R^d$ and the question is to study their stability and, when unstable,  how the dynamics can generate a moving Dirac mass.  

The Dirac steady states are given by
$$
\bar n(x;y)= \underline \rho(y) \delta (x-y).
$$
The total population size $ \underline \rho(y)$ is defined in both models by the constraint
$$
R \big(y, \underline I(y) \big) =0,
$$
with respectively for \eqref{paranm}  and \eqref{localcnm} 
$$
\underline I(y))= \psi(y) \underline \rho(y), \qquad  \text{resp. } \;  \underline I(y)=   \underline \rho(y) C(y,y) .
$$
A monotonicity assumption in $I$ for model \eqref{paranm}, namely $R_I(x,I)<0$ shows uniqueness of $\underline I(y)$ for a given $y$. In case of \eqref{localcnm} it is necessary that $r(y)>0$ for the positivity of $ \underline \rho(y) $.    
\\

In both models a 'strong' perturbation in measures  is stable, i.e. only on the weight; for $n^0=\rho^0 \delta (x-y)$, the solution is obviously 
$n(t,x)= \rho(t) \delta (x-y)$ with
$$
\f d{dt}  \rho(t)=  \rho(t) R\big(y, \psi(y)  \rho(t)\big),  \qquad  \text{resp. } \; \f d{dt}  \rho(t)=  \rho(t) [r(y)- \rho(t) C(y,y)],
$$
and 
$$
\rho(t) \underset{t \to \infty}{\longrightarrow}  \;  \underline \rho(y).
$$
This simple remark explains why, giving $I$ now,  the hypersurface $\{x, \; R(x,I)=0\}$  is a natural candidate for the location of a possible Dirac curve (see the introduction). 
\\

Apart from this stable one dimensional manifold, the Dirac steady states are usually unstable by perturbation in the weak topology. A way to quantify this instability is to follow the lines of \cite{OD.PJ.SM.BP:05} and consider at $t=0$ an exponentially concentrated initial data
$$
n^0(x)= e^{u^0_\ep(x)/\ep} \underset{\ep \to 0}{\longrightarrow}   \;  \underline \rho( \bar x^0) \delta (x-  \bar x^0).
$$
It is convenient to restrict our attention to $u^0_\ep$ uniformly concave, having in mind the gaussian case mentioned in the introduction. 
Then, $\ep$ measures the deviation from the  initial Dirac state and to see motion it is necessary to consider long times as $t/\ep$ or, equivalently, rescale the equation as  
$$
\ep \; \p_t n_\ep=  {n_\ep} R \big(x, I_\ep(t,x) \big),
$$
and our goal is to prove that 
$$
n_\ep(t,x) \underset{\ep \to 0}{\longrightarrow}   \bar n(t,x)= \bar \rho(t) \delta \big(x- \bar x(t) \big).
$$

Also the deviation to a Dirac state turns out to stay at the same size for all times  and we can better analyze this phenomena  using the WKB ansatz \eqref{wkb}. 
Indeed, $u_\ep$ satisfies the equation 
$$
\p_t u_\ep=  R \big(x, I_\ep(t,x) \big).
$$
As used by \cite{LD.PJ.SM.GR:08}, because $u^0$ is concave, assuming $x \mapsto R(x, I_\ep(t,x))$ is also concave (this only relies on assumptions on the data), we conclude that $u_\ep(t,\cdot)$ is also concave and thus has a {\em unique} maximum point $\bar x_\ep(t)$. The Laplace formula shows that, with $\bar x(t)$ the strong limit of $\bar x_\ep(t)$, 
$$
\f{n_\ep(t,x)}{\int_{\R^d}n_\ep(t,x) dx}  \underset{\ep \to 0}{\longrightarrow}     \delta \big(x- \bar x(t) \big).
$$

With some functional analysis, we are able to pass to the strong limit in $I_\e$ and $u_\ep$. Despite its nonlinearity, we find the same limiting equation,
\begin{equation} \label{hjnm}
\p_t u =  R \big(x, I(t,x) \big), \qquad u(t=0)= u^0.
\end{equation}
Still following the idea introduced in \cite{OD.PJ.SM.BP:05}, we may see $I(t,x)$ or $\rho(t)$ as a Lagrange multiplier for the constraint
\begin{equation} \label{constraintnm}
\max_{\R^d} u(t,x)= 0 = u\big(t, \bar x(t)\big),
\end{equation}
which follows from the  a priori bound $0<\rho(t) \leq \rho_M< \infty$. The mathematical justification of these developments is rather easy here. For    model \eqref{paranm} it uses a $BV$ estimate proved in \cite{GB.BP:07}. For model \eqref{localcnm} one has to justify persistence (that is $\rho_\ep$ stays uniformly positive) and strong convergence of $\rho_\ep(t)$. All this work is detailed below with the additional Laplace terms.

The constraint \eqref{constraintnm} allows us to recover the 'fast' dynamics of $I(t)$ and $\rho(t)$. Indeed, combined with  \eqref{hjnm}, it yields
\begin{equation} \label{levelsetnm}
R \big( \bar x(t), I(t) \big)=0,    \qquad  \text{resp. } \;  R \big( \bar x(t), I(t,  \bar x(t)) \big)=0.
\end{equation}
Assuming regularity on the data, $u(t,x)$ is three times differentiable and the constraint \eqref{constraintnm}  also gives
$$
\nabla u\big(t, \bar x(t)\big) =0. 
$$
Differentiating in $t$, we establish the analogue of the canonical equations in \cite{UD.RL:96} (see also  \cite{NC.RF.GB:01, OD.PJ.SM.BP:05, Raoul2009, SM.lisbon})
\begin{equation} \label{cenm}
\dot{ \bar x }(t)= \left( -D^2u \big(t, \bar x(t)\big) \right)^{-1}. D_x R\big(\bar x(t),\bar I(t, \bar x(t)) \big), \quad \bar x(0)= \bar x^0,
\end{equation}
(in the case of model \eqref{localcnm}, this means the derivative with respect to $x$ in both places). Inverting $I$ from the identities \eqref{levelsetnm} gives an autonomous equation. 
\\

This differential equation has a kind of Lyapunov functional and this makes it easy to analyze its long time behaviour. It is closely related
 to know what are the stable Dirac states for the weak topology; the so-called Evolutionary attractor or Convergence Stable Strategy in the language of adaptive dynamics \cite{OD:04, Raoul2009}. Eventhough this is less visible, it also carries regularity on the Lagrange multiplier $\rho(t)$ which helps for the functional analytic work in the case with mutations. Another use of \eqref{cenm} is to explain why, generically, only one pointwise Dirac mass can be sustained; as we explained earlier, the equation on $\dot{ \bar x }(t)$ also gives the global unknown $I(t)$ by coupling with \eqref{levelsetnm} and this constraint is very strong. See section \ref{sec:numerics} for an example.
\\

To conclude this quick presentation, we notice that the time scale (in $\ep$ here) has to be precisely adapted to the specific initial state under consideration. The initial state itself also has to be 'exponentially' concentrated along with our construction; this is the only way to observe the regular motion of the Dirac concentration point. This is certainly implicitly used in several works where such a behaviour is displayed, at least numerically. Of course, there are many other ways to concentrate the initial state with a 'tail' covering the full space so as to allow that  any trait $x$ can emerge; these are not covered by the present analysis.

\section{Competition through a single resource: assumptions and main results }
\label{sec:gca}

As used by \cite{LD.PJ.SM.GR:08}, concavity assumptions on the function $u_\ep$ in  \eqref{wkb} are enough to ensure concentration of $n_\ep(t, \cdot)$ as a {\em single} Dirac mass. We follow this line and make the necessary assumptions.

We start with assumptions on $\psi$:
\begin{equation} \label{aspsi}
0<\psi_m\leq \psi \leq \psi_M < \infty, \qquad \psi \in W^{2,\infty}(\R^d) .
\end{equation}

The assumptions on $R\in C^2$ are that  there is a constant $ I_M>0$  such that (fixing the origin in $x$ appropriately)
\begin{equation} \label{asrmax}
\max_{x \in \R^d} R(x,I_M) = 0 = R(0,I_M) ,
\end{equation}
\begin{equation} \label{asr}
-\underline{K}_1 |x|^2 \leq R(x,I) \leq \overline{K}_0 -\overline{K}_1 |x|^2,  \qquad \text{for }\;  0\leq I \leq I_M ,
\end{equation}
\begin{equation} \label{asrD2}
- 2\underline{K}_1 \leq D^2 R(x,I) \leq - 2\overline{K}_1 < 0    \text{ as symmetric matrices for }\, 0\leq I \leq I_M,
\end{equation}
\begin{equation} \label{asrDi}
- \underline{K}_2\leq \dis{\frac{\p R}{\p I} \leq - \overline{K}_2}, \qquad  \Delta (\psi R) \geq  - K_3.
\end{equation}
At some point we will also need that (uniformly in $0\leq I \leq I_M$)
\begin{equation} \label{asrD3}
D^3R(\cdot,I) \in  L^\infty(\R^d).
\end{equation}

\

Next the initial data $n_\ep^0$ has to be chosen compatible with the assumptions on $R$ and $\psi$. We require that there is a constant $I^0$ such that
\begin{equation} \label{asI}
0<I^0\leq I_\ep(0) := \int_{\R^d} \psi(x) n_\ep^0(x) dx< I_M,
\end{equation}
that we can write
$$
n_\ep^0=e^{u_\ep^0/\ep},  \qquad \text{with }\;  u_\ep^0 \in C^2(\R^d) \quad (\text{uniformly in } \ep),
$$
and we assume uniform concavity on $u^\ep$ too. Namely, there are positive constants $\underline{L}_0, \overline{L}_0,\underline{L}_1, \overline{L}_1$ such that
\begin{equation} \label{asu}
-\underline{L}_0 -\underline{L}_1 |x|^2 \leq u_\ep^0(x) \leq  \overline{L}_0 - \overline{L}_1 |x|^2 ,
\end{equation}
\begin{equation} \label{asuD2}
-2\underline{L}_1 \leq D^2u_\ep^0  \leq - 2\overline{L}_1 .
\end{equation}
We also make the following assumption which states that  the population is not initially maladapted and guarantees its survival
\begin{equation}
\label{as:survival}
\left( \frac{1}{\ep}\int \psi(x) R(x,I_\e(0)) n_\ep^0(x)\,dx \right)_- =o(1), \qquad \text{as $\e\to 0$}.
\end{equation}
For  Theorems \ref{th.cano} and \ref{th.long}  we also need that
\begin{equation} \label{asuD3}
D^3u_\ep^0 \in L^\infty(\R^d) \quad \text{componentwise uniformly in }\ep ,
\end{equation}
\begin{equation} \label{asuIni}
n_\ep^0(x) \underset{\ep\to 0}{\longrightarrow}  \bar{\rho}^0 \; \delta \big(x-\bar x^0\big) \text{ weakly in the sense of measures.}
\end{equation}

\

Next we need  to restrict the class of initial data to fit with $R$ through some compatibility conditions
\begin{equation} \label{asru}
4 \overline{L}_1^2 \leq \overline{K}_1 \leq  \underline{K}_1  \leq 4 \underline{L}_1^2.
\end{equation}

\

In the concavity framework of these assumptions, we are going to prove the following
\begin{theorem}[Convergence]\label{th.conv}  Assume \re{aspsi}-\re{asrDi}, \re{asI}-\re{as:survival} and \re{asru}. Then for all $T>0$, there is a $\ep_0>0$ such that for $\ep < \ep_0$ and $t \in [0,T]$, the solution $n_\ep$ to \eqref{para} satisfies,
\begin{align} \label{h}
&0< \rho_m \leq \rho_\ep(t):=\int_{\R^d} n_\ep\,dx \leq \rho_M+C\ep^2,\qquad 0< I_m \leq I_\ep(t) \leq I_M+C\ep^2 \quad \text{a.e.}
\end{align}
for some constant $\rho_m$, $I_m$. Moreover, $I_\e$ is uniformly bounded in $BV(\R^+)$ and after extraction of a subsequence $I_\ep$
\begin{equation} \label{limI}
I_\ep(t) \underset{\ep\to 0}{\longrightarrow}  \bar I(t)\quad \text{ in } L_{loc}^1(\R^+), \quad I_m \leq \bar I(t) \leq I_M \quad \text{a.e.},
\end{equation}
and
$ \bar I(t) $ is non-decreasing.
We also have weakly in the sense of measures for a subsequence $n_\ep$
\begin{equation} \label{limN}
n_\ep(t,x) \underset{\ep\to 0}{\longrightarrow}  \bar \rho(t) \; \delta \big(x-\bar x (t)\big).
\end{equation}
Finally, the pair $\big(\bar x(t),  \bar I(t) \big)$ also satisfies
\begin{equation}\label{eq.R0}
R\big(\bar x(t),\bar I(t) \big)=0 \quad \text{a.e.}
\end{equation}
\end{theorem}

In particular, there can be an initial layer on $I_\ep$ that makes a possible rapid variation of $I_\ep$ at $t\approx 0$ so that the limit satisfies $R(\bar x^0, I(0^+))=0$, a relation that might not hold true, even with $O(\ep)$, at the level of $n_\ep$.

\begin{theorem}[Form of canonical equation]\label{th.cano}
Assume \re{aspsi}-\re{asru}. Then,  $\bar x (\cdot)$ belongs to \\
$W^{1,\infty}(\R^+; \R^d)$ and satisfies
\begin{align}\label{eq.cano}
\dot{ \bar x }(t)= \left( -D^2u \big(t, \bar x(t)\big) \right)^{-1}. \nabla_x R\big(\bar x(t),\bar I(t) \big), \quad \bar x(0)=\bar{x}^0,
\end{align}
with $u(t,x)$ a $C^2$-function given below in \re{eq.u}, $D^3u \in L^\infty(\R^d)$, and initial data $\bar{x}^0$ given in \re{asuIni}.
Furthermore, we have $\bar I(t) \in W^{1,\infty}(\R^+)$.
\end{theorem}

We insist that the Lipschitz continuity at $t=0$ is with the value $I(0)= \lim_{t \to 0^+} I(t) \neq \lim_{\ep \to 0} I_\ep^0$; the equality might hold if the initial data is well-prepared.

\begin{theorem}[Long-time behaviour]\label{th.long}  With the assumptions \re{aspsi}-\re{asru}, equation \eqref{eq.cano} has a kind of Lyapunov functional, the limit $\bar I(t)$ is increasing and
\begin{equation} \label{limt}
\bar I(t) \underset{t\to \infty}{\longrightarrow}   I_M, \quad  \bar{x}(t) \underset{t\to \infty}{\longrightarrow}   \bar{x}_\infty=0.
\end{equation}
Finally, the limit is identified by $\na R(\bar{x}_\infty=0, I_M)=0$ (according to \eqref{asrmax}).
\end{theorem}

It is an open question to know if the full sequence converges. This is to say if the solution to the Hamilton-Jacobi equation is unique. The only uniqueness case in \cite{GB.BP:07} assumes a very particular form of $R(\cdot, \cdot)$.

\section{A-priori bounds on $\rho_\ep$, $I_\ep$ and their limits}
\label{se.rho}

Here, we establish the first statements of Theorem \ref{th.conv}.
As in \cite{GB.SM.BP:09} we can show with \re{asrmax} and \re{asrDi} that $I_\ep\leq I_M + C \ep^2$.
With \re{aspsi} (the bounds on $\psi$), we also have that
\begin{align}\label{eq.rho}
\rho_\ep(t) \leq I_M/ \psi_m + C \ep^2.
\end{align}

To achieve the lower bound away from $0$ is more difficult. We multiply the equation \re{para} by $\psi$ and integrate over $\R^d$, to arrive at
\begin{align}\label{eq.I_ep}
\ddt I_\ep(t) = \frac{1}{\ep}\int \psi R n_\ep\,dx + \ep\int n_\ep\De \psi  \,dx.
\end{align}
We define $\dis{J_\ep (t):=\frac{1}{\ep}\int \psi R n_\ep\,dx}$ and calculate its time derivative
\begin{align*}
\ddt J_\ep(t) = \frac{1}{\ep}\int \psi R \left(\frac{1}{\ep} R n_\ep +  \ep\De n_\ep\right)\,dx+ \frac{1}{\ep} \int \psi n_\ep\frac{\p R}{\p I}\left(J_\ep+ \ep \int n_\ep\Delta\psi \,dy \right)\,dx.
\end{align*}
So we estimate it from below using \re{aspsi}, \re{asrDi} and \re{eq.rho} by
\begin{align*}
\ddt J_\ep(t) \geq -C + \frac{1}{\ep}J_\ep \left(t \right) \int \psi n_\ep\frac{\p R}{\p I}\,dx,
\end{align*}
and we may bound the negative part of $J_\ep$ by
\begin{align}\label{eq.J}
\ddt (J_\ep(t))_- \leq C - \frac{\overline{K}_2}{\ep} I_\ep(t)  (J_\ep(t))_- .
\end{align}
Now for $\ep$ small enough, we can estimate $I_\ep(t)$ as
\begin{align}\label{eq.I_below}
I_\ep(t) &= I_\ep(0) + \int_0^t \dot{I_\ep}(s)\,ds = I_\ep(0) + \int_0^t J_\ep(s)\,ds + O(\ep)
\geq I^0/2 - \int_0^t  (J_\ep(s))_- \,ds,
\end{align}
and plugging this in the estimate \eqref{eq.J} leads to
\begin{align*}
\ddt (J_\ep(t))_- \leq C - \frac{\overline{K}_2}{\ep}\left(  I^0/2 - \int_0^t  (J_\ep(s))_- \,ds \right)  (J_\ep(t))_- .
\end{align*}
Now for $T>0$ fixed. If there exists $T'\leq T$ such that  $\int_0^{T'}  (J_\ep(s))_- \,ds =  I^0/4$, then we have
\begin{align*}
\ddt (J_\ep(t))_- \leq C - \frac{\overline{K}_2}{\ep}\f{ I^0}{4}  (J_\ep(t))_-  , \qquad 0\leq t\leq T'.
\end{align*}
Thus we obtain
\begin{align*}
(J_\ep(t))_- \leq (J_\ep(t=0))_- e^{-\overline{K}_2I^0t/(4\ep)} + \frac{4C\ep}{\overline{K}_2I^0} \left(1-e^{-\overline{K}_2I^0t/(4\ep)}  \right).
\end{align*}
Then, thanks to Assumption \re{as:survival} and for $\ep <\ep_0(T)$ small enough, we conclude that such a $T'$ does not exist {\it i.e.}
\begin{align}\label{eq.bvfirst}
\int_0^{T}  (J_\ep(s))_- \,ds \leq  I^0/4.
\end{align}

So from \re{eq.I_below}, we obtain
\begin{align}\label{eq.I_below1}
I_\ep(t) & \geq I^0/4, \qquad (J_\ep(t))_-  \underset{\ep\to 0}{\longrightarrow} 0 \qquad \text{ a. e. in }[0,T].
\end{align}
This also gives the lower bound $\rho_m \leq \rho_\ep(t)$ with $\rho_m := I^0/(4\psi_M)$.
\\

Finally, the estimate \eqref{eq.bvfirst} and the $L^\infty$ bounds on $I_\ep(t)$ give us a local BV bound, which will eventually allow us to extract a convergent subsequence for which \re{limI} holds. The obtained limit function $\bar I(t)$ is non-decreasing because in the limit the right-hand side of \re{eq.I_ep} is almost everywhere non-negative.

%
\section{Estimates on $u_\ep$ and its limit $u$}

In this section we introduce the major ingredient in our study, the function $u_\ep:= \ep \ln(n_\ep)$. We calculate
\begin{align*}
\p_t n_\ep= n_\ep \p_t  u_\ep/\ep, \quad \na n_\ep= n_\ep \na  u_\ep/\ep,  \quad \Delta n_\ep= n_\ep \Delta  u_\ep/\ep +n_\ep |\na  u_\ep|^2/\ep^2.
\end{align*}
Plugging this in \re{para}, we obtain that $u_\ep$ satisfies the Hamilton-Jacobi equation

\begin{equation}\label{eq.u_ep}
\left\{
\begin{array}{l}
\p_t u_\ep = |\na u_\ep|^2 + R(x,I_\ep(t)) + \ep \Delta u_\ep, \qquad x \in \R^d, \; t \geq 0,
\\[2mm]
u_\ep(t=0)= \ep \ln(n^0_\ep):=u_\ep^0.
\end{array}
\right.
\end{equation}

Our study of the concentration effect relies mainly on the asymptotic analysis of the family $u_\ep$ and in particular on its uniform regularity. We will pass to the (classical) limit in \eqref{eq.u_ep}, and this relies on  the

\begin{lemma}\label{le.u}
With the assumptions of Theorem \ref{th.conv}, we have  for $t\geq 0$,
\begin{align}  \label{eq.le.u}
-\underline{L}_0 -\underline{L}_1 |x|^2  -\ep 2d\underline{L}_1t \leq  u_\ep(t,x) \leq \overline{L}_0 -\overline{L}_1 |x|^2 +\left( \overline{K}_0 + 2 d \ep \overline{L}_1\right)t ,
\end{align}
\begin{align}  \label{eq.le.D2u}
-2\underline{L}_1 \leq D^2u_\ep(t,x)  \leq -2\overline{L}_1 .
\end{align}
\end{lemma}

This Lemma relies on a welknown (and widely used) fact that the Hamilton-Jacobi equations have a regime of regular solutions with concavity assumptions, \cite{GB:94,Lionsbook:82}.

\subsection{Quadratic estimates on $u_\ep$}
\label{se.quad}

First we achieve an upper bound,
defining  $\overline{u}_\ep(t,x):= \overline{L}_0 -\overline{L}_1 |x|^2 + C_0(\ep)t$ with $C_0(\ep) := \overline{K}_0 + 2 d \ep \overline{L}_1$, we obtain thanks to  \re{asr}, \re{asu} and \re{asru} that $\overline{u}_\ep(t=0) \geq u_\ep^0$ and
\begin{align*}
\p_t \overline{u}_\ep - |\na \overline{u}_\ep|^2 - R(x,I_\ep) - \ep \Delta \overline{u}_\ep \geq C_0(\ep) - 4 \overline{L}_1^2 |x|^2 - \overline{K}_0 +\overline{K}_1 |x|^2 - 2 d \ep \overline{L}_1\geq 0  .
\end{align*}
Next for the lower bound, we
define $\underline{u}_\ep(t,x):= -\underline{L}_0 -\underline{L}_1 |x|^2  -\ep C_1t $ with $C_1 := 2d\underline{L}_1$, we have

\noindent
$\underline{u}_\ep(t=0) \leq u_\ep^0$ and
\begin{align*}
\p_t \underline{u}_\ep - |\na \underline{u}_\ep|^2 - R(x,I_\ep) - \ep \Delta \underline{u}_\ep \leq -\ep C_1 -4\underline{L}_1^2 |x|^2 + \underline{K}_1 |x|^2 + \ep 2d\underline{L}_1 \leq 0.
\end{align*}

This concludes the proof of the first part of Lemma \ref{le.u} {\it i.e.} inequality \ref{eq.le.u}.

%
\subsection{Bounds on $D^2u_\ep$}
\label{se.reg}

We show that the semi-convexity and the concavity of the initial data is preserved by equation \re{eq.u_ep}.
For a unit vector $\xi$, we use the notation $u_\xi := \nabla_{\xi} u_\ep$ and $u_{\xi\xi} := \nabla^2_{\xi\xi} u_\ep$ to obtain
\begin{align*}
&u_{\xi t} =  R_\xi(x,I) +2\na u \cdot \na u_\xi   + \ep \Delta u_\xi,\\
&u_{\xi\xi t} = R_{\xi\xi}(x,I) + 2\na u_\xi \cdot \na u_\xi+2\na u \cdot \na u_{\xi\xi}   +  \ep \Delta u_{\xi\xi}.
\end{align*}
The first step is to obtain a lower bound on the second derivative {\it i.e.} semi-convexity.
It can be obtained in the same way as in \cite{GB.BP:08}:
Using  $|\na  u_\xi| \geq |u_{\xi\xi}|$ and the  definition $\underline{w}(t,x) :=  \min_\xi u_{\xi\xi}(t,x)$
leads to the inequality
\begin{align*}
\p_t \underline{w} \geq -2\underline{K}_1+ 2\underline{w}^2 +2\na u \cdot \na \underline{w}  +  \ep \Delta \underline{w}.
\end{align*}
By  a comparison principle and assumptions \re{asuD2}, \re{asru}, we obtain
\begin{align} 
\underline{w}\geq -2\underline{L}_1.
\end{align}

At every point $(t,x) \in \R^+\times\R^d$, we can choose an orthonormal basis such that $D^2u_\ep(t,x)$ is diagonal because it is a symmetric matrix. 
So we can estimate the mixed second derivatives in terms of $u_{\xi\xi}$. In particular, for each element $\xi$ of the latter basis,
we have $\na u_\xi=u_{\xi\xi} \xi$ and $|\na u_\xi|= |u_{\xi\xi}|$.

This enables us to show  concavity in the  next step.
We  start from the definition $\overline{w}(t,x) :=  \max_\xi u_{\xi\xi}(t,x)$ and the  inequality
\begin{align*}
\p_t \overline{w} \leq -2\overline{K}_1+ 2\overline{w}^2 +2\na u \cdot \na \overline{w}  +  \ep \Delta \overline{w}.
\end{align*}
By  a comparison principle and assumptions \re{asuD2}, \re{asru}, we obtain
\begin{align} 
\overline{w}\leq  -2\overline{L}_1.
\end{align}

From the space regularity gained and \re{asuD2}, we obtain $\na u_\ep$ locally uniformly bounded and thus from \re{eq.u_ep} for $\ep <\ep_0 $ that $\p_t u_\ep$ is locally uniformly bounded.

\subsection{Passing to the limit}
\label{se.limit}

From the regularity obtained in section \ref{se.reg}, it follows that we can extract a subsequence such that, for all $T>0$,
\begin{align*}
u_\ep(t,x)\underset{\ep\to 0}{\longrightarrow}u(t,x) \text{ strongly in } L^\infty \left(0,T;W_{loc}^{1,\infty}(\R^d)\right),
\end{align*}
\begin{align*}
u_\ep(t,x) \underset{\ep\to 0}{\xrightharpoonup{\quad} } u(t,x) \text{ weakly-* in } L^\infty\left(0,T;W_{loc}^{2,\infty}(\R^d)\right)\cap W^{1,\infty}\left(0,T;L_{loc}^{\infty}(\R^d)\right),
\end{align*}
and
\begin{equation} \label{eq:concavity_u}
-\underline{L}_0 -\underline{L}_1 |x|^2 \leq u(t,x) \leq  \overline{L}_0 - \overline{L}_1 |x|^2 + \overline{K}_0t,\quad
-2\underline{L}_1 \leq D^2u(t,x) \leq -2\overline{L}_1 \quad \text{ a.e.}
\end{equation}
\begin{equation} \label{eq:lipschitz_u}
u\in W^{1,\infty}_{\rm loc} (\R^+ \times \R^d).
\end{equation}
Notice that the uniform $W_{loc}^{2,\infty}(\R^d)$ regularity also allows to differentiate the equation in time, and find
$$
\f{\p^2}{\p t^2} u = \f{\p}{\p I} R\big(x,I(t)\big) \f{dI(t)}{dt}+ 2 \na u. [\na R(x,I(t))+ D^2u .\na u] .
$$
This is not enough to have $C^1$ regularity on $u$.

We also obtain that $u$ satisfies in the viscosity sense (modified as in \cite{GB.BP:07,GB.BP:08}) the equation
\begin{equation}\label{eq.u}
\left\{
\begin{array}{l}
\f{\p}{\p t} u = R\big(x,I(t)\big) + |\na u|^2,
\\[4mm]
\max_{\R^d} u(t,x)=0.
\end{array}
\right.
\end{equation}
The  constraint that the maximum vanishes is achieved, as in \cite{GB.SM.BP:09}, from the a priori bounds on $I_\ep$ and \eqref{eq:concavity_u}.

In particular $u$ is strictly concave, therefore it has exactly one maximum. This proves \re{limN} {\it i.e.} $n$ stays monomorphic and characterizes the Dirac location by
\begin{equation}\label{eq.diraclocation}
\max_{\R^d} u(t,x)=0= u\big(t, \bar x (t) \big).
\end{equation}
Moreover, as in \cite{GB.BP:08}  we can achieve \re{eq.R0} at each Lebesgue point of $I(t)$.

This completes the proof of Theorem \ref{th.conv}.

\section{Canonical equation, time asymptotic}
\label{se.proLi}

With the additional assumptions \re{asrD3} and \re{asuD3}, we can write our form of the canonical equation and show  long-time behavior. To do so, we first show that the third derivative is bounded. This allows us to establish rigorously the canonical equation  while it was only formally given in \cite{OD.PJ.SM.BP:05, GB.BP:08}. From this equation, we obtain regularity on $\bar x$ and $\bar I$. For  long-time behavior we show that $\bar I$ is strictly increasing as long as $\na R\big(\bar x,\bar I\big)\neq 0$ and
 this is based on a kind of Lyapunov functional.

\subsection{Bounds on third derivatives of $u_\ep$}
\label{se.thder}

For the unit vectors $\xi$ and $\eta$, we use the notation $u_\xi := \nabla_{\xi} u_\ep$, $u_{\xi\eta} := \nabla^2_{\xi\eta} u_\ep$ and $u_{\xi\xi\eta} := \nabla^3_{\xi\xi\eta} u_\ep$ to derive
\begin{align*}
\p_t u_{\xi\xi\eta} = 4 \na u_{\xi\eta}\cdot \na u_\xi + 2 \na u_\eta \cdot \na u_{\xi\xi} + 2 \na u \cdot \na u_{\xi\xi\eta}+R_{\xi\xi\eta} + \ep \Delta u_{\xi\xi\eta}.
\end{align*}
Again we can fix a point $(t,x)$ and choose an orthonormal basis such that $D^2\left(\na_\eta u_\ep(t,x)\right)$ is diagonal.
Let us define
\begin{align*}
M_1(t) := \max_{x,\xi,\eta} u_{\xi\xi\eta} (t,x).
\end{align*}
Since $-u_{\xi\xi\eta} (t,x) = \na_{-\eta} u_{\xi\xi} (t,x)$, we have $\dis{M_1(t) = \max_{x,\xi,\eta} \left|u_{\xi\xi\eta} (t,x)\right|}$. So with the maximum principle we obtain
\begin{align*}
\f{d}{dt} M_1 \leq 4dM_1\|D^2u_\ep\|_\infty +2dM_1\|D^2u_\ep\|_\infty  + R_{\xi\xi\eta}.
\end{align*}
Assumption \re{asuD3} gives us a bound on $M_1(t=0)$. So we obtain a $L^\infty$-bound on the third derivative uniform in $\ep$.

\subsection{Maximum points of  $u_\ep$}
\label{se.regMax}

Now we wish to establish the canonical equation. We denote the maximum point of $u_\ep(t,\cdot)$ by $\bar x_\ep(t)$.

Since $u_\ep\in C^2$, at maximum points we have $\na u_\ep(t,\bar{x}_\ep(t))=0$ and thus
\begin{align*}
\ddt \na u_\ep \big(t,\bar{x}_\ep(t)\big)=0.
\end{align*}
The chain rule gives
\begin{align*}
\frac{\p}{\p t} \na u_\ep\big(t,\bar{x}_\ep(t)\big)+ D_x^2 u_\ep\big(t,\bar{x}_\ep(t)\big)\dot{\bar{x}}_\ep(t)=0,
\end{align*}
and using the equation \re{eq.u_ep}, it follows further that, for almost every $t$,
\begin{align*}
D_x^2 u_\ep\big(t,\bar{x}_\ep(t)\big)\dot{\bar{x}}_\ep(t)=-\frac{\p}{\p t} \na u_\ep\big(t,\bar{x}_\ep(t)\big)=-\na_x R\big(\bar{x}_\ep(t),I_\ep(t)\big) - \ep \Delta \na_x u_\ep.
\end{align*}
Due to the uniform in $\ep$ bound on $D^3u_\ep$ and $R\in C^2$, we can pass to the limit in this equation and arrive at
\begin{align*}
\dot{ \bar x }(t)= \left( -D^2u \big(\bar x(t), t \big) \right)^{-1}. \nabla_x R\big(\bar x(t),\bar I(t) \big) \quad a.e.
\end{align*}

But we have further regularity in the limit and obtain the equations in the classical sense. We first notice that,
from $R\big(\bar x(t),\bar I(t) \big)=0$ and the assumptions \re{asr}, $\bar x (t)$  is bounded in $L^{\infty}(\R^+)$.
So we obtain that  $\bar x (t)$  is bounded in $W^{1,\infty}(\R^+)$. Because $I \mapsto R(\cdot,I)$ is invertible, it follows that  $\bar I (t)$  is bounded in $W^{1,\infty}(\R^+)$; more precisely we may differentiate the relation \re{eq.R0} (because $R\in C^2$) and find a differential equation on $I(t)$ that will be used later:
$$
\dot{ \bar x }(t) \cdot \na_x R + \dot{ \bar I }(t)\na_I R =0 \quad a.e.
$$

This completes the proof of Theorem \ref{th.cano}.

%
\subsection{Long-time behaviour}
\label{se.long}

It remains to prove the long time behaviour stated in Theorem \ref{th.long}.

We start from the canonical equation
$$
\f d {dt} \bar x(t)= (-D^2 u)^{-1} \nabla R\big(\bar x(t),\bar I(t)\big),
$$
and use some kind of Lyapunov functional. We calculate
\begin{align*}
\f d {dt}R\big(\bar x(t),\bar I(t)\big) &= \nabla R\big(\bar x(t),\bar I(t)\big)  \f d {dt} \bar x(t)+ R_I \big(\bar x(t),\bar I(t)\big)  \f {d\bar I} {dt} \\
&=\nabla R\big(\bar x(t),\bar I(t)\big)(-D^2 u)^{-1} \nabla R\big(\bar x(t),\bar I(t)\big)  + R_I \big(\bar x(t),\bar I(t)\big)  \f {d \bar I} {dt} .
\end{align*}

Now we also know from \re{eq.R0} that the left hand side vanishes. Then, we obtain
$$
\f {d} {d  t}\bar I(t)= \f{-1}{R_I \big(\bar x(t),\bar I(t)\big)}\nabla R\big(\bar x(t),\bar I(t)\big)(-D^2 u)^{-1} \nabla R\big(\bar x(t),\bar I(t)\big)  \geq 0.
$$
The inequality is strict as long as $\bar I(t) < I_M$.
Consequently, we recover that $\bar  I(t)$ is non-decreasing, as $t \to \infty$, $\bar I(t)$ converges, and subsequences of $\bar x(t)$ converge also (recall that $\bar x (t)$ is bounded). But we discover that  the only possible limits are such that $\nabla R(\bar x_\infty, I_\infty)= 0$. With relation \re{eq.R0}, assumptions \re{asrmax} and \re{asrDi} this identifies the limit as announced in Theorem \ref{th.long}.

\section{Numerics}
\label{sec:numerics}

\graphicspath{{./Images/}}
\begin{figure}[h]
\centerline{  \subfloat[asymmetric I.D.]{\includegraphics[width=0.5\textwidth]{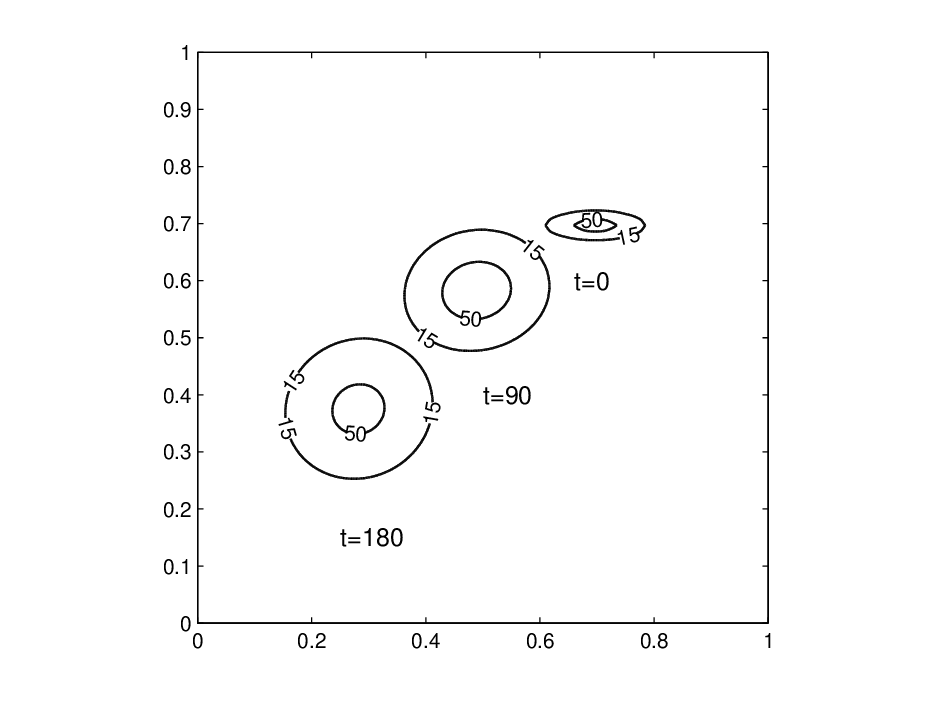}} \qquad
\subfloat[symmetric I.D.]{\includegraphics[width=0.5\textwidth]{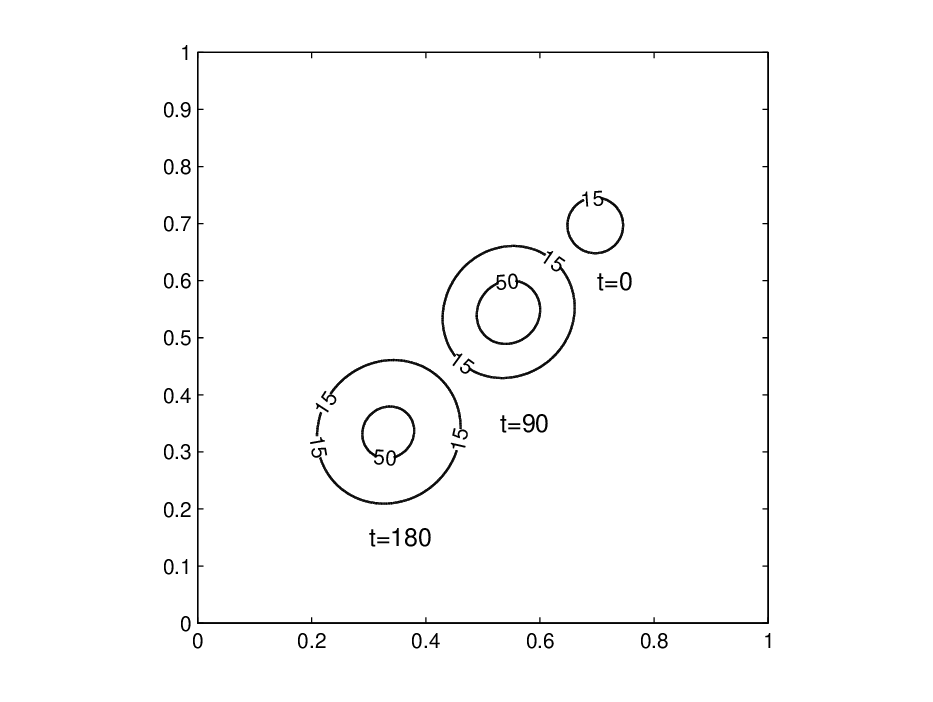}}}
\caption{Dynamics of the density $n$ with asymmetric  initial data \re{eq.iniasy} (left) and symmetric initial data \re{eq.inisy} (right).
These computations illustrate the effect of the matrix $(-D^2u)^{-1}$ in the dynamics of the concentration point according to the form of canonical equation \re{eq.cano}. The plots show the level sets $\lbrace (x,y) | n(x,y)=15 \rbrace$ and $\lbrace (x,y) | n(x,y)=50\rbrace$.
}
\label{fi.asymIni}
\end{figure}

The canonical equation is not self-contained because the effect of mutations appears through the matrix $(-D^2u)^{-1}$.  Nevertheless it can be used to explain several effects. The purpose of this section is firstly to illustrate how it acts on the dynamics, secondly to see the effect of $\ep$ being not exactly zero, and thirdly to explain why it is generic, in high dimensions as well as in one dimension \cite{GB.BP:08} that pointwise Dirac masses (and not on curves) can exist.
\\

We first illustrate the fact that a isotropic approximation of a Dirac mass will give rise to  different dynamics than an anisotropic. This anisotropy is measured with $u$ and we choose two initial data. In the first case $-D^2u^0$ is ''far'' from the identity matrix and in the second case it is isotropic:
\begin{align}\label{eq.iniasy}
n_0^\ep(x,y)=C_{\text{mass}} \exp(-(x-.7)^2/\ep -12(y-.7)^2/\ep),
\end{align}
\begin{align}\label{eq.inisy}
n_0^\ep(x,y)=C_{\text{mass}} \exp(-2.4(x-.7)^2/\ep -2.4(y-.7)^2/\ep).
\end{align}

We also choose a growth rate $R$ with gradient along the diagonal:
\begin{align}
R(x,y,I)=2- I- 0.6(x^2+y^2).
\end{align}
Here although, we start with initial data centered on the diagonal and $\na R$ pointing along the diagonal to the origin, the concentration point with the anisotropic initial data \eqref{eq.iniasy} leaves the diagonal  (cf. Figure \ref{fi.asymIni} (a)).  The isotropic initial data moves along the diagonal as expected by symmetry reasons (cf. Figure \ref{fi.asymIni} (b)).

The numerics has been performed in Matlab with parameters as follows. 
The plots show the level sets $\lbrace (x,y) | n(x,y)=15 \rbrace$ and $\lbrace (x,y) | n(x,y)=50\rbrace$
corresponding to $t$ = 0, 90 and 180 (time in units of $dt$): $\ep$ is chosen to be $0.005$ and $C_{\text{mass}}$ such that the initial mass in the computational domain is equal to $0.3$. The equation is solved  by an implicit-explicit finite-difference method on square grid consisting of $100 \times 100$ points  (time step $dt= 0.005$).
\\

\begin{figure}[h]
\centering
\includegraphics[width=0.45\textwidth]{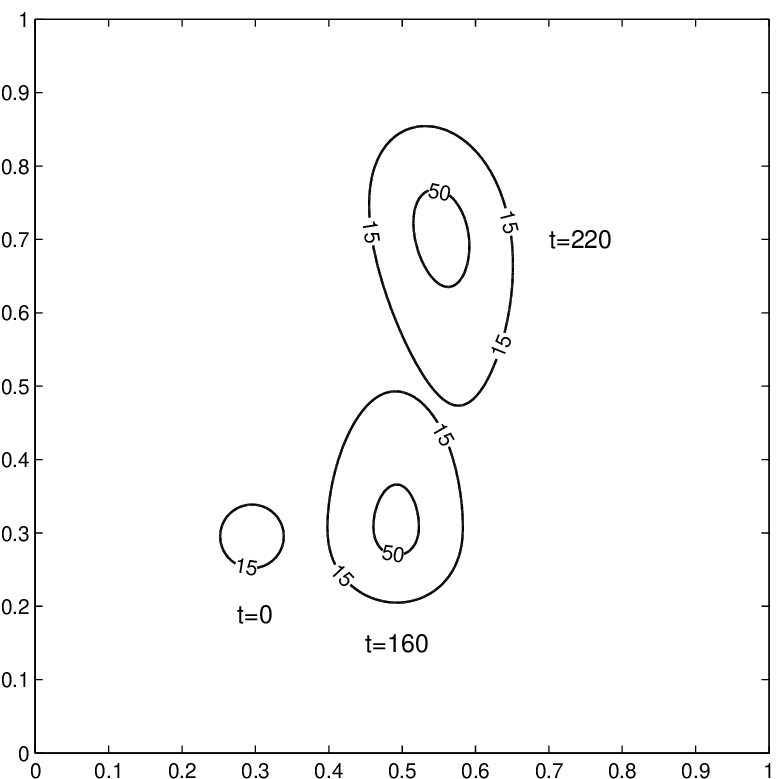}
\caption{
This figure illustrates the effect of $\ep$ being not exactly zero. The dynamics of the density $n$ with symmetric initial data is plotted for $t$ = 0, 160 and 220 in units of $dt$ and the limiting behavior is a motion along the axis $y=0.3$. The plot shows the level sets $\lbrace (x,y) | n(x,y)=15 \rbrace$ and $\lbrace (x,y) | n(x,y)=50\rbrace$.
}
\label{fi.playR}
\end{figure}

The second example is to illustrate the role of the parameter $\ep$ for symmetric initial data:
\begin{align}
n_0^\ep(x,y)=C_{\text{mass}} \exp(-(x-.3)^2/\ep -(y-.3)^2/\ep),
\end{align}
\begin{align}
R(x,y,I(t))=0.9-I+5(y-.3)_+^2+2.3(x-.3)\quad \text{ with } I(t):= \int n(t,x)\,dx.
\end{align}
In this example, we start with symmetric initial data centered on the line $y = 0.3$ and the gradient of $R$ along this line ($y=0.3$) is $(1,0)$. Hence, the canonical equation in the limit $\ep =0$ predicts a motion in the $x$ direction on this line. One however observes in Figure \ref{fi.playR} that the maximum point leaves this line because $\ep$ does not vanish. Notice that $\p_y R\equiv 0$ below the line $y=0.3$.

In this computation, performed with Matlab,  $\ep$ is chosen to be $0.004$, $C_{\text{mass}}$ such that the initial mass in the computational domain is equal to $0.3$ and square grid consisting of $150 \times 150$ points (time step $dt= 8.8889\cdot10^{-4}$).
\\
\graphicspath{{./Images/}}

\begin{figure}\centerline{
\subfloat[asymmetric R]{\includegraphics[width=0.4\textwidth]{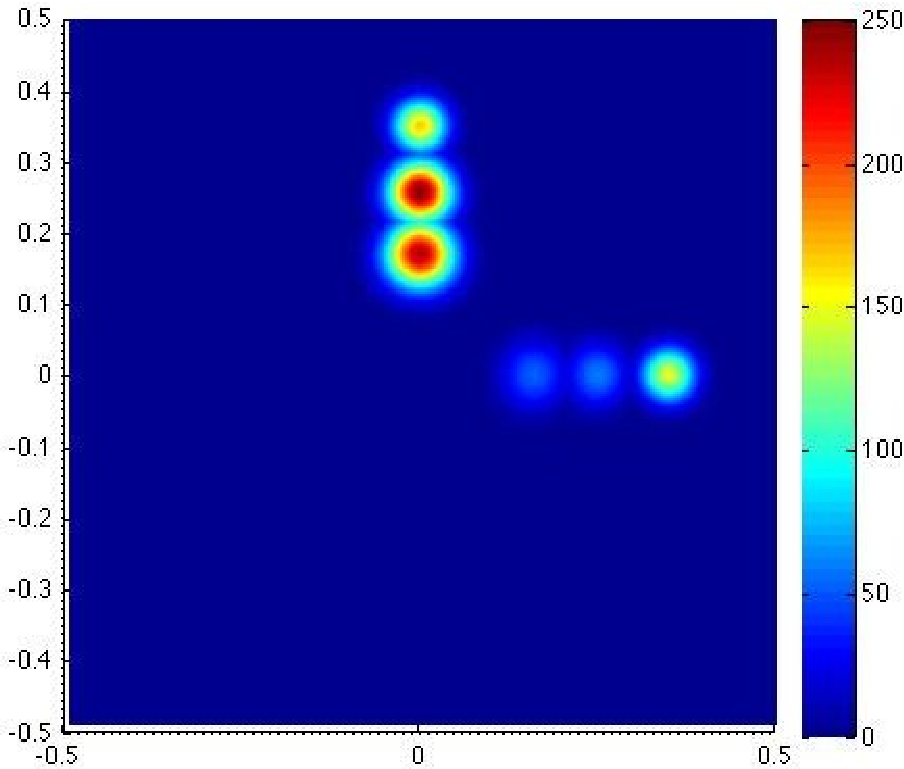}} \qquad
\subfloat[symmetric R]{\includegraphics[width=0.4\textwidth]{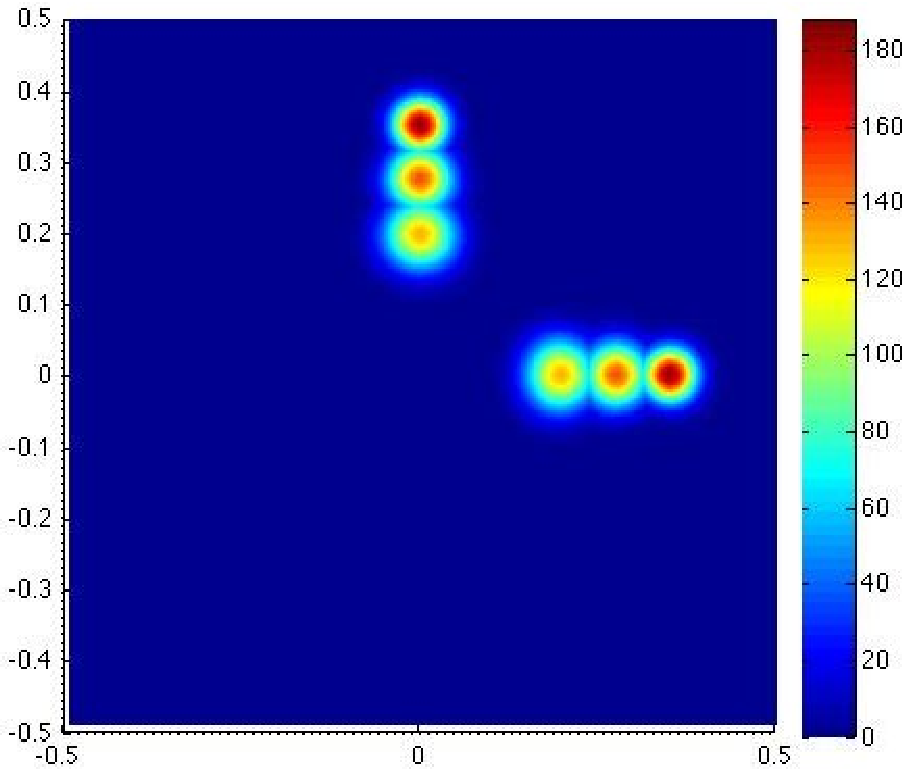}}  }
\caption{ This figure illustrates that, except for particular symmetry conditions, a single Dirac mass is exhibited by Lotka-Volterra equations. We depict the density $n$ with asymmetric (left) and symmetric (right) growth rate $R$ plotted for $t$ = 5, 90 and 180 in units of $dt$.}
\label{fi.R_asy_sy}
\end{figure}

With our third example we wish to illustrate that, except in particular symmetric geometries, only a single Dirac mass can 
be sustained by the Lotka-Volterra equations with a single resource in the framework of section \ref{sec:gca}. We place initially two symmetric deltas on the $x$ and the $y$-axis:
\begin{align}
n_0^\ep(x,y)=C_{\text{mass}} \left[ \exp\left(-\frac{2.4}{\ep}\left((x-.25\sqrt{2})^2+ y^2\right)\right)+\exp\left(-\frac{2.4}{\ep}\left((y-.25\sqrt{2})^2 +x^2\right)\right)\right],
\end{align}
We seek for asymmetry in the growth rate $R$ under the form
\begin{align}
R(x,y, I)=3-1.5 I-5.6(y^2+R_e x^2).
\end{align}
In the special case $R_e=1$, all isolines of $R$ are circles then the two concentration points just move symmetrically to the origin cf. Figure \ref{fi.R_asy_sy} (b).
However, if we choose $R_e=1.1$ {\it i.e.} all isolines of $R$ are ellipses then one of  the two concentration points disappears cf. Figure \ref{fi.R_asy_sy} (a).
The intuition behind is that the canonical equation \re{eq.cano} should hold for the two points. However the  constraint \re{eq.R0} given by $\rho(t)$ is the same for the two points and this is a contradiction. One of the two points has to disappear right away.

The numerics is performed with  $\ep = 0.003$ and $C_{\text{mass}}$ such that the initial mass in the computational domain is equal to $0.3$. The equation is solved  by an implicit-explicit finite-difference method on square grid consisting of $100 \times 100$ points (time step $dt= 0.001$).

\graphicspath{{./Images/}}

%
\section{Extension: non-constant diffusion}
\label{sec:extension}

Our results can be extended to include a diffusion coefficient depending on $x$. This leads to the equation
\begin{equation} \label{parab}
\p_t n_\ep - \ep \na  \cdot \left( b(x) \na n_\ep\right)=  \frac{n_\ep}{\ep} R \big(x, I_\ep(t) \big),  \qquad t>0, \; x \in \R^d.
\end{equation}
Our assumptions on $b$ are that there are positive constants $b_m$, $b_M$, $B_1$, $B_2$ and $B_3$ such that
\begin{align}\label{asb}
b_m\leq b\leq b_M, \qquad |\na b(x)| \leq B_1 \frac{1}{1 + |x|}, \qquad \left|Tr(D^2b(x))\right| \leq B_2 \frac{1}{(1 + |x|)^2}, \qquad  \left|D^3b\right| \leq B_3.
\end{align}

Our assumptions on the initial data and on $R$ are the same, as before \re{aspsi}--\re{asuIni}. However, we have to supplement the assumption \re{asrDi} to take $b$ into account:
\begin{equation} \label{asrDi1}
\na b \cdot \na (\psi R)  \geq  - K_3.
\end{equation}
These assumptions will in the following allow us to obtain a gradient bound
\begin{align}\label{eq.u_epgrad}
|\na u_\ep(t,x) | \leq C_{\na u} (1 + |x|).
\end{align}
This bound enables us to formulate the compatibility conditions which replace \re{asu}: we need
\begin{align}\label{asK}
B_2C_{\na u}^2-2\overline{K}_1 <0
\end{align}
and define
\begin{align*} 
\overline{K}_b:=\frac{2B_1-\sqrt{4B_1^2 - 2 b_M\left(B_2C_{\na u}^2-2\overline{K}_1\right)}}{b_M},
\end{align*}
\begin{align*} 
\underline{K}_b:=\frac{-2B_1-\sqrt{4B_1^2 + 2 b_m\left(B_2C_{\na u}^2+2\underline{K}_1\right)}}{b_m},
\end{align*}
to require
\begin{equation} \label{asuD2b}
-\underline{K}_b \leq D^2u_\ep^0  \leq -\overline{K}_b,
\end{equation}
\begin{equation} \label{asrub}
4 b_M \overline{L}_1^2 \leq \overline{K}_1 \leq  \underline{K}_1  \leq 4 b_m \underline{L}_1^2.
\end{equation}

Our goal is to prove the following

\begin{theorem}[Convergence]\label{th.convb}  Assume \re{aspsi}-\re{asrDi}, \re{asu}, \re{as:survival}, \re{asb}, \re{asrDi1}, \re{asK}, \re{asuD2b} and \re{asrub}. Then the solution $n_\ep$ to \eqref{parab} satisfies
for all $T>0$, for  $\ep$ small enough and $t \in [0,T]$
\begin{align} \label{hb}
&0< \rho_m \leq  \rho_\ep(t) \leq \rho_M+C\ep^2,\qquad I_m \leq I_\ep(t) \leq I_M+C\ep^2 \quad \text{a.e.}
\end{align}
Moreover, there is a subsequence $I_\ep$ such that
\begin{equation} \label{limIb}
I_\ep(t) \underset{\ep\to 0}{\longrightarrow}  \bar I(t)\quad \text{ in } L_{loc}^1(\R^+), \quad I_m \leq \bar I(t) \leq I_M \quad \text{a.e.},
\end{equation}
and
$ \bar I(t) $ is non-decreasing.
Furthermore, we have weakly in the sense of measures for a subsequence $n_\ep$
\begin{equation} \label{limNb}
n_\ep(t,x) \underset{\ep\to 0}{\longrightarrow}  \bar \rho(t) \; \delta \big(x-\bar x (t)\big),
\end{equation}
and
the pair $\big(\bar x(t),  \bar I(t) \big)$ also satisfies
\begin{equation}\label{eq.R0b}
R\big(\bar x(t),\bar I(t) \big)=0 \quad \text{a.e.}
\end{equation}

\end{theorem}

\begin{theorem}[Form of canonical equation]\label{th.canob}

With the assumptions \re{aspsi}-\re{asu}, \re{as:survival}--\re{asuIni},   \re{asb}, \re{asrDi1}, \re{asK}, \re{asuD2b} and \re{asrub}, $\bar x$ is a $W^{1,\infty}(\R^+)$-function satisfying
\begin{align}\label{eq.canob}
\dot{ \bar x }(t)= \left( -D^2u \big(t,\bar x(t) \big) \right)^{-1}. \nabla_x R\big(\bar x(t),\bar I(t) \big), \quad \bar x(0)=\bar{x}^0, 
\end{align}
with $u(t,x)$ a $C^2$-function given below in \re{eq.ub}, $D^3u \in L^\infty(\R^d)$, and initial data $\bar{x}^0$ given in \re{asuIni}.
Furthermore, we have $\bar I(t) \in W^{1,\infty}(\R^+)$.
\end{theorem}
The end of this section is devoted to the proof of these Theorems.
The a priori bounds \re{hb}, \re{limIb} on $\rho_\ep$ and $I_\ep$ can be established as before.

As before we study the function $u_\ep:= \ep \ln(n_\ep)$. We obtain that $u_\ep$ satisfies the Hamilton-Jacobi equation
\begin{equation}\label{eq.u_epb}
\left\{
\begin{array}{l}
\p_t u_\ep =   R \big(x, I_\ep(t) \big) +b|\na u_\ep|^2+\ep \na b \cdot \na u_\ep+ \ep b\Delta u_\ep,  \qquad t>0, \; x \in \R^d,\\
u_\ep(t=0)= \ep \ln(n^0_\ep).
\end{array}
\right.
\end{equation}

In order to adapt our method to this equation we need a bound on the gradient of $u_\ep$. We achieve this following arguments in  \cite{PL:85, GB.SM.BP:09}:

\paragraph{gradient bound}
Let us define $v(t,x)$ by $u_\ep= K_v - v^2$ where we choose $K_v$  large enough to have $v >\delta >0$ on $[0,T]$ uniform in $\ep$.  We calculate
\begin{align*}
\na u_\ep = -2 v \na v \quad \text{ and } \quad\De u_\ep = -2 v \De v -2 |\na v|^2
\end{align*}
and obtain from \re{eq.u_epb}
\begin{align}\label{eq.v}
-2v \p_t v = R + 4b v^2 |\na v|^2 -2 \ep \na b \cdot v \na v - 2\ep v b \De v - 2 \ep b |\na v|^2  .
\end{align}
Dividing by $-2v$, taking the derivative with respect to $x_i$ and defining $p := \na v$, we have
\begin{multline*}
\p_t p_i = -\left(\frac{R}{2v}\right)_{x_i}-2 p_i b |p|^2 - 2 v b_{x_i} |p|^2 - 4vb p \cdot \na p_i + \ep \na b_{x_i} \cdot p + \ep \na b \cdot \na p_i \\+ \ep b_{x_i} \De v + \ep b \De p_i + \ep \frac{b_{x_i} }{v} |p|^2 -  \ep \frac{b }{v^2} |p|^2 p_i +  2\ep \frac{b }{v} p \cdot \na p_i.
\end{multline*}
Multiplying \re{eq.v} by $\dis{\frac{b_{x_i}}{2bv}}$ and adding to the equation above, we obtain
\begin{multline*}
\p_t \left(p_i - \frac{b_{x_i}v}{b}\right) = -\left(\frac{R}{2v}\right)_{x_i}-2 p_i b |p|^2 - 2 v b_{x_i} |p|^2 - 4vb p \cdot \na p_i + \ep \na b_{x_i} \cdot p + \ep \na b \cdot \na p_i \\ + \ep b \De p_i + \ep \frac{b_{x_i} }{v} |p|^2 -  \ep \frac{b }{v^2} |p|^2 p_i +  2\ep \frac{b }{v} p \cdot \na p_i + \frac{b_{x_i}R}{2bv} +2b_{x_i}v|p|^2-\ep\frac{b_{x_i}}{b}  \na b \cdot p-\ep\frac{b_{x_i}}{v}|p|^2.
\end{multline*}
Now we define
\begin{align}
M_b(t):= \max_{i,x}\left[(p_i)_-, (p_i)_+\right] \geq 0.
\end{align}
If $\max_{i,x}(p_i)_- \leq \max_{i,x}(p_i)_+$, we have
\begin{multline*}
\p_t \left(M_b - \frac{b_{x_i}v}{b}\right) \leq C-2 b_m M_b^3 +  2 |v| |b_{x_i}| d^2 M_b^2  + \ep d \left| \na b_{x_i}\right| M_b  \\  + \ep \frac{|b_{x_i} |}{\delta} d^2M_b^2 +C +2|b_{x_i}||v| d^2 M_b^2+\ep\frac{|b_{x_i}|}{b_m}d |\na b| M_b  +\ep\frac{|b_{x_i}|}{\delta}d^2 M_b^2.
\end{multline*}
Since $\dis{\frac{b_{x_i}v}{b}}$ is bounded, we have $M_b$ bounded from above.

If $\max_{i,x}(p_i)_- > \max_{i,x}(p_i)_+$, we show similarly a bound on $M_b$ and therefore achieve \re{eq.u_epgrad}.

To prove the concavity and semi-convexity results, we only give formal arguments for the limit case. To adapt the argument for the $\ep$-case is purely technical:

For a  unit vector $\xi$, we define $u_{\xi}:= \na_\xi u_\ep$ and $u_{\xi\xi}: =\na_{\xi\xi} u_\ep$ to obtain
\begin{align}\label{eq.ub_xi}
\p_t u_\xi = R_{\xi}+b_\xi |\na u|^2 + 2 b \na u \cdot \na u_\xi ,
\end{align}
and
\begin{align}\label{eq.ub_xixi}
\p_t u_{\xi\xi} = R_{\xi\xi}+b_{\xi\xi} |\na u|^2 + 4 b_\xi \na u \cdot \na u_\xi+2 b \na u \cdot \na u_{\xi\xi} +2 b | \na u_{\xi}|^2 .
\end{align}
With the definition  $\overline{w}(t,x) :=  \max_\xi u_{\xi\xi}(t,x)$ and assumptions \re{asb} we have
\begin{align*} 
\p_t \overline{w} \leq -2\overline{K}_1+B_2C_{\na u}^2 + 4 B_1 C_{\nabla u}|\overline{w}| +2 b \na u \cdot \na \overline{w} +2 b_M \overline{w}^2.
\end{align*}
With assumption \re{asK}, $0$ is a supersolution to
\begin{align*} 
\p_t \overline{w}^* = -2\overline{K}_1+B_2C_{\na u}^2 - 4 B_1  C_{\nabla u}\overline{w}^* +2 b \na u \cdot \na \overline{w}^* +2 b_M (\overline{w}^*)^2,
\end{align*}
so we know from assumption \re{asuD2b} that $\overline{w}  \leq 0$.
Therefore it follows further that
\begin{align*} 
\overline{w} \leq \overline{K}_b 
.
\end{align*}

For the lower bound, we use the definition  $\underline{w}(t,x) :=  \min_\xi u_{\xi\xi}(t,x)$ and the inequality
\begin{align*} 
\p_t \underline{w} \geq -2\underline{K}_1-B_2C_{\na u}^2 - 4 B_1  C_{\nabla u} |\underline{w}| +2 b \na u \cdot \na \underline{w} +2 b_m \underline{w}^2.
\end{align*}
Since we already know that  $\underline{w} \leq 0$, we obtain
\begin{align*} 
\underline{w} \geq \underline{K}_b.
\end{align*}
We can achieve this at the $\ep$-level using the equation
\begin{align}\label{eq.ub_xiep}
\p_t u_\xi = R_{\xi}+b_\xi |\na u|^2 + 2 b \na u \cdot \na u_\xi +\ep \na b_\xi \cdot \na u + \ep \na b \cdot \na u_\xi + \ep b_\xi \De u + \ep b \Delta u_\xi,
\end{align}
and
\begin{multline}\label{eq.ub_xixiep}
\p_t u_{\xi\xi} = R_{\xi\xi}+b_{\xi\xi} |\na u|^2 + 4 b_\xi \na u \cdot \na u_\xi+2 b \na u \cdot \na u_{\xi\xi} +2 b | \na u_{\xi}|^2 \\+\ep\na b_{\xi\xi}\cdot \na u + 2 \ep\na b_\xi\cdot \na u_\xi +\ep\na b \cdot \na u_{\xi\xi}+ \ep b_{\xi\xi} \De u+ 2\ep b_\xi \De u_\xi + \ep b \De u_{\xi\xi}.
\end{multline}
Now we define
\begin{align*}
f:=\frac{2b_\xi}{b} \text{ and } g:=\frac{bb_{\xi\xi}-2b_\xi^2}{b^2},
\end{align*}
multiply \re{eq.ub_xiep} by $f$, substract it from \re{eq.ub_xixiep}, multiply \re{eq.u_epb} by $g$, substract it  
to obtain
\begin{multline} 
\p_t (u_{\xi\xi}-fu_\xi-gu) = R_{\xi\xi} +2 b \na u \cdot \na u_{\xi\xi} +2 b | \na u_{\xi}|^2 +\ep \na b_{\xi\xi}\cdot \na u\\ + 2 \ep \na b_\xi \na u_\xi + \ep \na b \cdot \na u_{\xi\xi} + \ep b \De u_{\xi\xi} - f R_{\xi}  -\ep f \na b_\xi \cdot \na u-\ep f \na b \cdot \na u_\xi-gR-\ep g\na b \cdot \na u.
\end{multline}
The remaining steps can be done similar as before. For the Hamilton-Jacobi-equation on $u$, we obtain the variant 
\begin{equation}\label{eq.ub}
\left\{
\begin{array}{l}
\p_t u =   R \big(x, \bar I(t) \big) +b(x) \; |\na u|^2,
\\[2mm]
\dis \max_{x \in \R^d} u(t,x)= 0, \qquad \forall t \geq 0.
\end{array}
\right.
\end{equation}

%
\section{Direct competition}
\label{sec:localcompetition}

The other class of models we handle are populations with direct competition kernel $C(x,y) \geq 0$, that is
\begin{align}\label{para'}
\p_t n_\ep(t,x) = \frac{1}{\ep}n_\ep(t,x) \left(r(x)- \int_{\R^d} C(x,y)n_\ep(t,y) \,dy\right) + \ep \De n_\ep(t,x).
\end{align}
The term $r(x)$ is the intra-specific growth rate (and has a priori no sign) and the integral term models an additional  contribution to the death rate due to competition between different traits. Notice that the choice $C(x,y)= \psi(y) \Phi(x)$ will reduce this model to a particular case of those in \eqref{para}. This class of model is also very standard, see \cite{RC.JH:05,GM.MG:05,NC.RF.SM:06,NC.RF.SM:08,Raoulphd} and the references therein.
We call it direct competition in opposition to more realistic models where competition is through resources \cite{SM.BP.JW:10}.
\\

For the initial data, we assume as before \re{asu}--\re{asuIni}. Concerning $r(x)$ and $C(x,y)$ we assume $C^1$ regularity and that there are constants $\rho_M>0$, $\underline{K}'_1>0$... such that

\begin{align}\label{ascbasic}
C(x,x) >0 \qquad   \forall x \in \R^d,
\end{align}
\begin{align}\label{asco}
\int_{\R^d}\int_{\R^d}n(x)C(x,y)n(y) \,dy dx \geq \frac{1}{\rho_M}  \int_{\R^d}n(x) \,dx \int_{\R^d}r(x)n(x) \,dx \qquad \forall n \in L^1_+(\R^d).
\end{align}
This assumption is weaker than easier conditions of the type
$$
C(x,y) \geq \frac{1}{\rho_M} r(x)  \quad \text{or } \; C(x,y) \geq \frac{1}{2 \rho_M} [r(x) +r(y)] .
$$
Because it is restricted to positive functions, it is a pointwise positivity condition on $C(x,y)$ in opposition to the positivity as operator that occurs for the entropy method in \cite{PJ.GR:09}. 

Then, we make again concavity assumptions. Namely that concavity on $r$ is strong enough to compensate for concavity in $C$
\begin{equation} \label{asr'}
-\underline{K}'_1 |x|^2 \leq r(x)-\sup_y C(x,y) \rho_{M}  \leq r(x) \leq \overline{K}'_0 -\overline{K}'_1 |x|^2,
\end{equation}
\begin{equation} \label{asrD2'}
- 2\underline{K}'_1 \leq D^2r(x)-\sup_y\big(D^2C(x,y)\big)_+\rho_{M}  \leq D^2r(x)+\sup_y\big(D^2C(x,y)\big)_-\rho_{M}  \leq - 2\overline{K}'_1 ,
\end{equation}
as symmetric matrices, where the positive and negative parts are taken componentwise. As for regularity, we will use
\begin{equation} \label{asrD3'}
D^3r-\sup_y\big(D^3C(\cdot,y)\big)_+\rho_{M}, \qquad D^3r +\sup_y\big(D^3C(\cdot,y)\big)_-\rho_{M}  \in  L^\infty(\R^d).
\end{equation}

The initial data  is still supposed to concentrate at a point $\bar x^0$ following  \re{asu}--\re{asuIni}. But because persistence, i.e. that $n_\e$ does not vanish, is more complicated  to control, we need two new conditions
\begin{equation} \label{rx0}
r(\bar x^0)> 0, 
\end{equation}
\begin{equation} \label{ascompuo}
\int_{\R^d}  n_\ep(t,x) dx \leq \rho^0_{M}.
\end{equation}
We also need a compatibility condition with $R$
\begin{equation} \label{asru'}
4 \overline{L}_1^2 \leq \overline{K}'_1 \leq  \underline{K}'_1  \leq 4 \underline{L}_1^2 .
\end{equation}

The interpretation of our assumptions is that the intra-specific growth rate $r$ dominates strongly the competition kernel. This avoids the branching patterns that are usual in this kind of models \cite{GM.MG:05, SG.VV.PA:06,BP.SG:07,Raoulphd}. Our concavity assumptions also implies that there is no continuous solution $N$ to the steady state equation without mutations  $N(x) \left(r(x)- \int_{\R^d} C(x,y) N(y) \,dy\right) =0$. This makes a difference with the entropy method used in \cite{PJ.GR:09} as well as the positivity condition on the kernel that, compared to \eqref{asco}, also involves $r(x)$.

Our goal is to prove the following results
\begin{theorem}[Convergence]\label{th.conv'}
With the assumptions \re{ascbasic}--\re{asrD2'} and  \re{asu}, \re{ascompuo}--\re{asru'},  the solution $n_\ep$ to \eqref{para'} satisfies,
\begin{align} \label{h'}
&0\leq \rho_\ep(t):= \int_{\R^d} n_\ep(t,x) dx \leq \rho_{M}\quad \text{a.e.}
\end{align}
and there is a subsequence  such that
\begin{equation} \label{limI'}
\rho_\ep(t) \underset{\ep\to 0}{\xrightharpoonup{\quad}}  \bar \rho(t)\quad \text{ in weak-}\star \;  L^\infty(\R^+), \qquad 0 \leq \bar \rho(t) \leq  \rho_{M}\quad \text{a.e.}
\end{equation}
Furthermore, we have weakly in the sense of measures for a subsequence $n_\ep$
\begin{equation} \label{limN'}
n_\ep(t,x) \underset{\ep\to 0}{\xrightharpoonup{\quad}}  \bar \rho(t) \; \delta \big(x-\bar x (t)\big), \qquad \f{n_\ep(t,x) }{\int_{\R^d} n_\ep(t,x)  dx} \underset{\ep\to 0}{\xrightharpoonup{\quad}} \delta \big(x-\bar x (t)\big),
\end{equation}
and
the pair $\big(\bar x(t),  \bar \rho(t) \big)$ also satisfies
\begin{align}\label{eq.R0'}
\bar \rho(t) \big[ r\big(\bar x(t)\big)-\bar \rho(t)C\big(\bar x(t),\bar x(t)\big) \big] \geq 0.
\end{align}
\end{theorem}

With the assumptions of Theorem \ref{th.conv'}, we do not know if $\rho_\ep$ converges strongly because we do not have the equivalent of the $BV$ quantity in Theorem \ref{th.conv}. We can only prove it with stronger assumptions. This is stated in the

\begin{theorem}[Form of canonical equation]\label{th.cano'}
We assume \re{asu}--\re{asuIni} and \re{asco}--\re{asru'}. Then, for the $C^2$-function $u(t,x)$ given below in \re{eq.u'ex} with $D_x^3u \in L^\infty_{\rm loc} \big(\R^+;  L^\infty(\R^d) \big)$, $\bar x \in W^{1,\infty}(\R^+)$  satisfies 
\begin{align}\label{eq.cano'}
\dot{\bar x}(t) = \left(-D^2u\big(t,\bar x(t)\big)\right)^{-1} \cdot \left[\na_x r\big(\bar x(t)\big) - \bar \rho(t) \na_x C\big(\bar x(t),\bar x(t)\big) \right],
\end{align}
with initial data $\bar{x}^0$ given in \re{asuIni}.
Furthermore, $\rho_\ep$ converges strongly and  we have $\bar \rho(t) \in W^{1,\infty}(\R^+)$, 
\begin{align}\label{eq.R0''}
r\big(\bar x(t)\big)-\bar \rho(t)C\big(\bar x(t),\bar x(t)\big)=0 ,
\end{align}
\begin{equation}\label{limrho}
r\big(\bar x(t)\big) \geq r(\bar x^0) e^{-Kt}, \qquad \bar \rho(t)  \geq \rho^0e^{-Kt}  .
\end{equation}
\end{theorem}

We may find some kind of gradient flow structure for the canonical equation when $C(x,y)$ is symmetric and obtain

\begin{theorem}[Long time behavior]\label{th.compltb}
We make the assumptions of Theorem \ref{th.cano'}, $C(x,y)=C(y,x)$ and 
\begin{align} \label{comp:logconcave}
x \mapsto  \Phi(x):= \ln r(x) - \ln C(x,x) \quad \text{ is strictly concave in the set } \; \{ r>0\} .
\end{align}
Then, as $t \to \infty$,  $\bar \rho(t) \to  \bar \rho_\infty >0$, $ \bar x(t) \to  \bar x_\infty$ and $ \bar x_\infty$ is the maximum point of $\Phi$.
\end{theorem}

\subsection{A-priori bounds on $\rho_\ep$}
\label{sec:comp_apb}

The main new difficulty with the competition model comes from a priori bounds on the total population. In particular it is not known if there are $BV$ quantities proving that $\rho_\ep(t)$ converges strongly. Even non extinction is not longer automatic.

One side of the inequality \re{h'} is given by  $n_\ep \geq 0$, for the other side we
integrate \re{para'} over $\R^d$ and use \eqref{asco} to write
\begin{align*}
\f{d}{d t} \int_{\R^d} n_\ep(t,x)\,dx &= \frac{1}{\ep} \int_{\R^d} n_\ep(t,x) r(x)\,dx- \frac{1}{\ep} \int_{\R^d}\int_{\R^d} n_\ep(t,x) C(x,y)n_\ep(t,y) \,dydx \\
& \leq \frac{1}{\ep} \int_{\R^d} n_\ep(t,x) r(x)\,dx\left[ 1 - \frac{\int_{\R^d} n_\ep(t,x)\,dx }{\rho_M}  \right],
\end{align*}
therefore (and even though $r$ can change  sign) we conclude thanks to \re{ascompuo}
\begin{align*}
\int_{\R^d} n_\ep(t,x)\,dx  \leq \rho_{M}.
\end{align*}

\subsection{Passing to the limit}

The proofs of the remaining parts of Theorems  are close to those already written before. We only give the main differences here. They rely again on the WKB ansatz $u_\ep:= \ep \ln(n_\ep)$. We obtain as before that $u_\ep$ satisfies the Hamilton-Jacobi equation
\begin{equation}\label{eq.u_ep'}
\left\{
\begin{array}{l}
\f{\p}{\p t} u_\ep (t,x) = r(x) -\dis{ \int_{\R^d}C(x,y)n_\ep(t,y) \,dy } + |\na u_\ep|^2  + \ep \De u_\ep,  \qquad t>0, \; x \in \R^d,\\
u_\ep(t=0)= u^0_\ep .
\end{array}
\right.
\end{equation}

Similarly to Lemma \ref{le.u} we can prove the

\begin{lemma}\label{le.u'}
With the assumptions of Theorem \ref{th.conv'}, we have for all $t\geq 0$
\begin{align*}
-\underline{L}_0 -\underline{L}_1 |x|^2  -\ep 2d\underline{L}_1t \leq  u_\ep(t,x) \leq \overline{L}_0 -\overline{L}_1 |x|^2 +\left( \overline{K}_0' + 2 d \ep \overline{L}_1\right)t ,
\end{align*}
\begin{align*}
-2\underline{L}_1 \leq D^2u_\ep(t,x)  \leq -2\overline{L}_1.
\end{align*}
\end{lemma}

\proof The first line holds because the right (resp. left) hand side of the inequality is a super (resp. sub) solution thanks to assumption \eqref{asr'} and using the control of $n_\e$ by $\rho_M$. For the second line, the upper and lower bound use the maximum principle on the equation for $D^2u_\e$ and the compatibility conditions \eqref{asru'} as in section \ref{se.reg}.
\qed

\

From the regularity obtained, it follows that we can extract a subsequence such that $u_\ep(t,x)\underset{\ep\to 0}{\longrightarrow}u(t,x)$, locally uniformly as in section \ref{se.limit}. We also obtain from Lemma \ref{le.u'}
\begin{equation}\label{ub'}
-\underline{L}_0 -\underline{L}_1 |x|^2 \leq u(t,x) \leq  \overline{L}_0 - \overline{L}_1 |x|^2+\overline{K}'_0t,\qquad
-2\underline{L}_1 \leq D^2u(t,x) \leq - 2\overline{L}_1\quad \text{ a.e.}
\end{equation}
and that $u$ satisfies, in the viscosity sense (modified as in \cite{GB.BP:07,GB.BP:08,GB.SM.BP:09}),  the equation
\begin{equation}\label{eq.u'}
\left\{
\begin{array}{l}
\f{\p}{\p t} u = r(x) -\bar{\rho}(t)C\big(x,\bar x(t)\big)+ |\na u|^2,
\\[2mm]
\max_{\R^d} u(t,x)\leq 0.
\end{array}
\right.
\end{equation}
The constraint is now relaxed to an inequality because we know that the total mass is bounded but we do not control the mass from below at this stage. In the framework of Theorem \ref{th.cano'}, we prove later on that the constraint is always an equality (see \eqref{eq.u'ex}).

It might be that $ \bar \rho(t)$ vanishes and then $\bar x (t)$ does not matter here, nevertheless we still have
\begin{equation}\label{eq.liconstraint}
\max_{\R^d} u(t,x) =  u\big(t, \bar x(t)\big).
\end{equation}

Using  the control \eqref{ub'} and this Hamilton-Jacobi equation  we obtain \re{limN'} with the same arguments as in Section \ref{se.limit}. The new difficulty is that $\rho(t)$ might vanish in particular when the constraint is strict  $\max_{\R^d} u(t,x) < 0$, an option that we will discard later. Because  of that, we also obtain the restriction on times in  \re{eq.R0'} which can be completed as (in the viscosity sense)
\begin{equation}\label{eq.R=0}
\f{d}{d t} u \big(t,\bar x(t)\big) = r\big(\bar x(t)\big)-\bar \rho(t)C\big(\bar x(t),\bar x(t)\big).
\end{equation}

We also have, as in section \ref{sec:comp_apb}, 
$$
\ep \; \f{d}{d t} \rho_\ep(t) = \int_{\R^d} n_\ep(t,x) r(x)\,dx- \int_{\R^d}\int_{\R^d} n_\ep(t,x) C(x,y)n_\ep(t,y) \,dydx.
$$
Passing to the weak limit (integration by parts and using boundedness of $\rho_\ep$),  we find that 
\begin{equation}\label{comp:weakl}
\bar \rho(t) r\big( \bar x (t)\big) =\text{ w-lim} \int_{\R^d}\int_{\R^d} n_\ep(t,x) C(x,y)n_\ep(t,y) \,dydx \geq \bar \rho(t)^2 C\big( \bar x (t),\bar x (t)\big).
\end{equation}
This proves \eqref{eq.R0'} and concludes the proof of Theorem \ref{th.conv'}.

\subsection{Form of the canonical equation}

We continue with the proof of Theorem \ref{th.cano'} and we begin with the derivation of \re{eq.cano'}.
\\

The third derivative of $u_\e$ is bounded using assumption \re{asrD3'} and following the same arguments in Section \ref{se.thder}. Then similarly to Section \ref{se.regMax}, we have the regularity $D^3_x u \in L^\infty\big((0,T)\times \R^d\big)$ for all $T>0$, $\f{\p}{\p t} u$ and $D^3_{txx} u \in L_{\rm loc}^\infty\big(\R^+\times \R^d\big)$.
\\

The canonical equation  \re{eq.cano'} is established a.e. as in section \ref{se.regMax} using the maximum points of $u_\ep$ and passing to the limit.  From  \re{h'}, \re{eq.cano'} and \re{ub'}, we next obtain that  $|\f{d}{dt}\bar x (t)|$  is uniformly bounded.

\subsection{Persistence}

Now we prove that $u\big(t,\bar x(t)\big)=0$ for all $t\geq 0$ and that $\bar \rho(t)>0$ a.e. $t$. We cannot obtain this directly and thus we begin with  proving $r(\bar x(t))> 0$. We prove indeed the first part of the inequality \re{limrho}.
\\

We prove this by contradiction. We suppose that $t_0$ is the first point such that $r(\bar x(t_0))=0$. We notice that $\bar x(t)$ being lipschitzian and using assumption \re{rx0}, for all $t<t_0$, we have $r(\bar x(t))>0$. Therefore with assumption \re{asr'} we deduce that $\bar x(t)$ remains bounded for $t\in [0,t_0]$. Using \re{eq.cano'} and \re{ub'} we have
\begin{align}\nonumber
\f {d}{dt}r(\bar x(t)) &=\nabla_x r(\bar x(t))\cdot \dot{\bar x}(t)\\\nonumber
&=\nabla_x r(\bar x(t))\cdot\left(-D^2u\big(t,\bar x(t)\big)\right)^{-1} \cdot \left[\na_x r\big(\bar x(t)\big) - \bar \rho(t) \na_x C\big(\bar x(t),\bar x(t)\big) \right]\\\nonumber
&\geq -\bar \rho(t) |\nabla_x r(\bar x(t))||\na_x C\big(\bar x(t),\bar x(t)\big)|.
\end{align}
Consequently, using \re{eq.R0'}, we obtain
$$\f {d}{dt}r(\bar x(t))\geq -r(\bar x(t)) |\nabla_x r(\bar x(t))|\f{|\na_x C\big(\bar x(t),\bar x(t)\big)|}{C\big(\bar x(t),\bar x(t)\big)}.$$
Moreover we know that $\bar x(t)$ remains bounded for $t\in [0,t_0]$ and thus we have \\
$\dis{\inf_{t\in [0,t_0]}C(\bar x(t), \bar x(t)) \geq \eta_2> 0}$.  We conclude that, for $K$ a positive constant, 
$$\f {d}{dt}r(\bar x(t)) \geq -Kr(\bar x(t)), \qquad \text{for }\;0\leq t\leq t_0.$$
Starting with $r(\bar x^0)>0$ according to \re{rx0}, this inequality is in contradiction with $r(\bar x(t_0))=0$. Therefore for all $t>0$ we have $r(\bar x(t))>0$ and thus this inequality is true for all $t>0$. Thereby we obtain the first part of \eqref{limrho}. From the latter and using again \fer{asr'}, we also deduce that $\bar x(t)$ remains bounded for all $t>0$. 
\\

Next we use \re{eq.R0'}, \eqref{eq.R=0} and the positivity of $r(\bar x(t))$ to obtain 
$$
u(t,\bar x(t)) -u^0(\bar x^0)=\int r(\bar x(t))-\bar \rho(t)C(\bar x(t),\bar x(t))dt \geq  \mathbf{1}_{\bar \rho(t)=0}r(\bar x(t)) \geq 0.
$$
We deduce, using \re{asuIni}, that $u(t,\bar x(t))=0$ for all $t\geq 0$. Thus the equation on $u$ is in fact
\begin{equation}\label{eq.u'ex}
\left\{
\begin{array}{l}
\f{\p}{\p t} u = r(x) -\bar{\rho}(t)C\big(x,\bar x(t)\big)+ |\na u|^2,
\\[2mm]
\max_{\R^d} u(t,x) = 0.
\end{array}
\right.
\end{equation}

This identity also proves the identity \eqref{eq.R0''} and  thus that  \eqref{comp:weakl} holds as an equality, which is equivalent to say that the weak limit of $\rho_\ep(t)$ is in fact a strong limit.
\\

We may now use \eqref{eq.R0''} to conclude that $\bar \rho (t)$  is also bounded in $W^{1,\infty}(\R^+)$. To do so, we first differentiate \re{eq.R0''} and find again some kind of gradient flow structure
\begin{align*}
\na r\big(\bar x(t)\big)\cdot\dot{\bar x}(t)-\dot{\bar \rho}(t)C\big(\bar x(t),\bar x(t)\big)-\bar \rho(t)\left[\na_x C\big(\bar x(t),\bar x(t)\big)+\na_y C\big(\bar x(t),\bar x(t)\big)\right]\dot{\bar x}(t)=0.
\end{align*}
With \re{eq.cano'}, it follows that
\begin{align}\label{comp:gradflow}
\dot{\bar x}(t)\cdot\left(-D^2u\right)\cdot \dot{\bar x}(t) = \dot{\bar \rho}(t)C\big(\bar x(t),\bar x(t)\big)+\bar \rho(t) \na_y C\big(\bar x(t),\bar x(t)\big)\dot{\bar x}(t).
\end{align}
From the uniform bounds proved before, there is a constant $\eta_3$ such that
\begin{align}   \label{eq.K'}
\frac{ \na_y C\big(\bar x(t),\bar x(t)\big)\cdot \dot{\bar x}(t)}{C\big(\bar x(t),\bar x(t)\big)} \leq \eta_3.
\end{align}
Using the latter and \re{ub'} we conclude that, for $K$ a positive constant, 
$$
\dot{\bar \rho}(t)\geq -K\bar \rho(t).
$$
Thus we obtain \re{limrho}. This completes the proof of Theorem \ref{th.cano'}.

\subsection{Long time behavior}

It remains to prove Theorem \ref{th.compltb}. Assuming that $C(x,y)$ is symmetric, we can find a quantity which is non-decreasing in time, which replaces the quantity $\bar I$ in section \ref{se.long}. We compute, from the relation \eqref{comp:gradflow}, 
\begin{align}   
\ddt \left[\bar{\rho}^2(t)C\big(\bar x(t),\bar x(t)\big)\right]\left(2\bar{\rho}(t)\right)^{-1}& = \dot{\bar \rho}(t) C\big(\bar x(t),\bar x(t)\big) + \bar{\rho}(t) (\na_y C)  \dot{\bar x}(t) \nonumber
\\
& = \dot{\bar x}(t)\cdot\left(-D^2u\right)\cdot \dot{\bar x}(t)  \geq 0.\label{comp:gradflow2}
\end{align}

As $t$ tends to infinity, we may consider a subsequence $t_n$ such that $\bar \rho(t_n) \to  \bar \rho_\infty$, $ \bar x(t_n) \to  \bar x_\infty$. From \eqref{comp:gradflow2}, we  may also assume  $\dot{ \bar x}(t_n) \to  0$. Therefore 
$$
\nabla r( \bar x_\infty) =   \bar \rho_\infty \nabla_x C(\bar x_\infty,\bar x_\infty), \qquad r( \bar x_\infty) =   \bar \rho_\infty C(\bar x_\infty,\bar x_\infty).
$$

From these relations, we first conclude that $\bar \rho_\infty$ is positive because $r$ is concave and its gradient only vanishes at a point where $r$ is positive.

Then we combine the relations and conclude that
$$
\f{\nabla r( \bar x_\infty)}{r( \bar x_\infty)} =\f{\nabla_x C(\bar x_\infty,\bar x_\infty)}{ C(\bar x_\infty,\bar x_\infty)}.
$$
The assumption \eqref{comp:logconcave} then concludes on the uniqueness of such a point $\bar x_\infty$ and thus on the convergence of $\bar x(t)$.
\qed

\section*{Acknowledgement}
This research is supported by Award No. KUK-I1-007-43, made by King
Abdullah University of Science and Technology (KAUST).

%
%


\begin{thebibliography}{10}

\bibitem{GB:94}
G.~Barles.
\newblock {\em Solutions de viscosit\'e des \'equations de
  {H}amilton-{J}acobi}, volume~17 of {\em Math\'ematiques \& Applications
  (Berlin) [Mathematics \& Applications]}.
\newblock Springer-Verlag, Paris, 1994.

\bibitem{GB.LE.PS:90}
G.~Barles, L.~C. Evans, and P.~E. Souganidis.
\newblock Wavefront propagation for reaction-diffusion systems of {PDE}.
\newblock {\em Duke Math. J.}, 61(3):835--858, 1990.

\bibitem{GB.SM.BP:09}
G.~Barles, S.~Mirrahimi, and B.~Perthame.
\newblock Concentration in {L}otka-{V}olterra parabolic or integral equations:
  a general convergence result.
\newblock {\em Methods Appl. Anal.}, 16(3):321--340, 2009.

\bibitem{GB.BP:07}
G.~Barles and B.~Perthame.
\newblock Concentrations and constrained {H}amilton-{J}acobi equations arising
  in adaptive dynamics.
\newblock {\em Contemp. Math.}, 439:57--68, 2007.

\bibitem{NC.RF.GB:01}
N.~Champagnat, R.~Ferri\`ere, and G.~Ben~Arous.
\newblock The canonical equation of adaptive dynamics: A mathematical view.
\newblock {\em Selection}, 2:73--83, 2001.

\bibitem{NC.RF.SM:06}
N.~Champagnat, R.~Ferri\`ere, and S.~M\'el\'eard.
\newblock Unifying evolutionary dynamics: From individual stochastic processes
  to macroscopic models.
\newblock {\em Theoretical Population Biology}, 69(3):297--321, 2006.

\bibitem{NC.RF.SM:08}
N.~Champagnat, R.~Ferri\`ere, and S.~M\'el\'eard.
\newblock {\em Individual-based probabilistic models of adaptive evolution and
  various scaling approximations}, volume~59 of {\em Progress in Probability}.
\newblock Birkha\"{u}ser, 2008.

\bibitem{NC.PJ:10}
N.~Champagnat and P.-E. Jabin.
\newblock The evolutionary limit for models of populations interacting
  competitively with many resources.
\newblock {\em Preprint}, 2010.

\bibitem{MC.CE.SM:19}
M.~Costa, C.~Etchegaray, and S.~Mirrahimi.
\newblock Survival criterion for a population subject to selection and
  mutations; application to temporally piecewise constant environments.
\newblock {\em Preprint}.

\bibitem{RC.JH:05}
R.~Cressman and J.~Hofbauer.
\newblock Measure dynamics on a one-dimensional continuous trait space:
  theoretical foundations for adaptive dynamics.
\newblock {\em Theoretical Population Biology}, 67(1):47--59, 2005.

\bibitem{LD.PJ.SM.GR:08}
L.~Desvillettes, P.-E. Jabin, S.~Mischler, and G.~Raoul.
\newblock On mutation-selection dynamics.
\newblock {\em Commun. Math. Sci.}, 6(3):729--747, 2008.

\bibitem{UD.RL:96}
U.~Dieckmann and R.~Law.
\newblock The dynamical theory of coevolution: {A} derivation from stochastic
  ecological processes.
\newblock {\em J. Math. Biol.}, 34:579--612, 1996.

\bibitem{OD:04}
O.~Diekmann.
\newblock A beginner's guide to adaptive dynamics.
\newblock In {\em Mathematical modelling of population dynamics}, volume~63 of
  {\em Banach Center Publ.}, pages 47--86. Polish Acad. Sci., Warsaw, 2004.

\bibitem{OD.PJ.SM.BP:05}
O.~Diekmann, P.-E. Jabin, S.~Mischler, and B.~Perthame.
\newblock The dynamics of adaptation: an illuminating example and a
  {H}amilton-{J}acobi approach.
\newblock {\em Th. Pop. Biol.}, 67(4):257--271, 2005.

\bibitem{LE:98}
L.~C. Evans.
\newblock {\em Partial differential equations}, volume~19 of {\em Graduate
  Studies in Mathematics}.
\newblock American Mathematical Society, Providence, RI, 1998.

\bibitem{LE.PS:89}
L.~C. Evans and P.~E. Souganidis.
\newblock A {PDE} approach to geometric optics for certain semilinear parabolic
  equations.
\newblock {\em Indiana Univ. Math. J.}, 38(1):141--172, 1989.

\bibitem{WF.PS:86}
W.~H. Fleming and P.~E. Souganidis.
\newblock {PDE}-viscosity solution approach to some problems of large
  deviations.
\newblock {\em Ann. Scuola Norm. Sup. Pisa Cl. Sci.}, 4:171--192, 1986.

\bibitem{SG.VV.PA:06}
S.~Genieys, V.~Volpert, and P.~Auger.
\newblock Pattern and waves for a model in population dynamics with nonlocal
  consumption of resources.
\newblock {\em Math. Model. Nat. Phenom.}, 1(1):63--80, 2006.

\bibitem{PJ.GR:09}
P.-E. Jabin and G.~Raoul.
\newblock Selection dynamics with competition.
\newblock {\em J. Math. Biol.}, To appear.

\bibitem{Lionsbook:82}
P.~L. Lions.
\newblock {\em Generalized solutions of {H}amilton-{J}acobi equations},
  volume~69 of {\em Research Notes in Mathematics}.
\newblock Pitman Advanced Publishing Program, Boston, 1982.

\bibitem{PL:85}
P.-L. Lions.
\newblock Regularizing effects for first-order {H}amilton-{J}acobi equations.
\newblock {\em Applicable Analysis}, 20:283--307, 1985.

\bibitem{SM.lisbon}
S.~M\'el\'eard.
\newblock Random modeling of adaptive dynamics and evolutionary branching.
\newblock In J.~F. Rodrigues and F.~Chalub, editors, {\em The Mathematics of
  Darwin's Legacy}, Mathematics and Biosciences in Interaction. Birkh\"auser
  Basel, 2010.

\bibitem{GM.MG:05}
G.~Mesz\'ena, M.~Gyllenberg, F.~J. Jacobs, and J.~A.~J. Metz.
\newblock Link between population dynamics and dynamics of darwinian evolution.
\newblock {\em Phys. Rev. Lett.}, 95(7):078105, Aug 2005.

\bibitem{M.P.B.M:10}
S.~Mirrahimi, B.~Perthame, E.~Bouin, and P.~Millien.
\newblock Population formulation of adaptative meso-evolution; theory and
  numerics.
\newblock In J.~F. Rodrigues and F.~Chalub, editors, {\em The Mathematics of
  Darwin's Legacy}, Mathematics and Biosciences in Interaction. Birkh\"auser
  Basel, 2010.

\bibitem{SM.BP.JW:10}
S.~Mirrahimi, B.~Perthame and J.~Y.~Wakano.
\newblock Evolution of species trait through resource competition.
\newblock  {\em Preprint}, 2010.

\bibitem{GB.BP:08}
B.~Perthame and G.~Barles.
\newblock Dirac concentrations in {L}otka-{V}olterra parabolic {PDE}s.
\newblock {\em Indiana Univ. Math. J.}, 57(7):3275--3301, 2008.

\bibitem{BP.SG:07}
B.~Perthame and S.~G{\'e}nieys.
\newblock Concentration in the nonlocal {F}isher equation: the
  {H}amilton-{J}acobi limit.
\newblock {\em Math. Model. Nat. Phenom.}, 2(4):135--151, 2007.

\bibitem{Raoulphd}
G.~Raoul.
\newblock {\em Etude qualitative et num\'{e}rique d'\'{e}quations aux
  d\'{e}riv\'{e}es partielles issues des sciences de la nature}.
\newblock PhD thesis, ENS Cachan, 2009.

\bibitem{Raoul2009}
G.~Raoul.
\newblock Local stability of evolutionary attractors for continuous structured
  populations.
\newblock {\em accepted in Monatsh. Math.}, 2010.

\end{thebibliography}
\end{document}